\documentclass[a4paper,12pt,dvipdfmx]{article}
\usepackage[dvips]{graphicx}
\usepackage{latexsym}
\usepackage{amssymb}
\usepackage{theorem}
\usepackage{amsmath}
\usepackage{amscd}
\usepackage{bm}
\usepackage{tikz}
\newtheorem{theorem}{Theorem}[section]
\newtheorem{corollary}[theorem]{Corollary}

\newtheorem{example}[theorem]{Example}
\newtheorem{proposition}[theorem]{Proposition}
\newtheorem{remark}[theorem]{Remark}
\newtheorem{definition}[theorem]{Definition}
\newcommand{\demo}{\par\noindent{\it Proof. \/}\ }
\newcommand{\enD}{\hfill $\Box$\vspace{3truemm} \par}
\newcommand{\R}{\mathbb{R}}

\newcommand{\rank}{\operatorname{rank}}

\newcommand{\bijar}[1][]{%
\ar[#1]
\ar@<0.7ex>@{}[#1]|-*[@]{\sim}} 
\begin{document}
\title{Circular evolutes and involutes of framed curves in the Euclidean space}
\author{Shun'ichi Honda and Masatomo Takahashi}
\date{\today}
\maketitle
%
%
\begin{abstract}
We introduce circular evolutes and involutes of framed curves in the Euclidean space. 
Circular evolutes of framed curves stem from the curvature circles of Bishop directions and singular value sets of normal surfaces of Bishop directions. 
On the other hand, involutes of framed curves are direct generalizations of involutes of regular space curves and frontals in the Euclidean plane. 
We investigate properties of normal surfaces, circular evolutes, and involutes of framed curves. 
We can observe that taking circular evolutes and involutes of framed curves are opposite operations under suitable assumptions, similarly to evolutes and involutes of fronts in the Euclidean plane. 
Furthermore, we investigate the relations among singularities of normal surfaces, circular evolutes, and involutes of framed curves. 
\end{abstract}
%
%
\renewcommand{\thefootnote}{\fnsymbol{footnote}}
\footnote[0]{2020 Mathematics Subject classification: 53A04, 53A05, 57R45, 58K05.}
\footnote[0]{Key Words and Phrases: Bishop frame, evolute, framed curve, framed surface, involute, normal surface, parallel curve, singular point.}
%
%
\section{Introduction}
For fronts in the Euclidean plane, the second author and Tomonori Fukunaga investigated that evolutes and involutes of fronts correspond to differential and integral of Legendre curvatures in classical calculus (cf. \cite{FT5}). 
On the other hand, the authors introduced spherical evolutes of framed curves in the Euclidean space that stem from osculating spheres and singular value sets of focal surfaces (cf. \cite{HT2}). 
The spherical evolutes of framed curves have certain geometric meanings but don't relate to involutes of framed curves. 
Involutes of framed curves considered here are direct generalizations of involutes of regular space curves and frontals in the Euclidean plane (cf. \cite{F1, FT4, FT5}). 
\par 
In this paper, we define circular evolutes of framed curves with respect to Bishop directions. 
They stem from curvature circles with respect to Bishop directions and singular value sets of normal surfaces with respect to Bishop directions. 
We investigate properties of normal surfaces, circular evolutes, and involutes of framed curves. 
\par 
Bishop directions are key role in order to define circular evolutes. 
The Frenet frame and Frenet-Serret formula are important for non-degenerate regular space curves and Frenet type framed base curves (cf. \cite{H1}). 
However, the Frenet frame is not a Bishop frame except for the torsion vanish, that is, space curves contain a plane. 
Actually, we cannot define circular evolutes by using Frenet frame except for the torsion vanish. 
Hence, the pair of a curve and frame is an important concept (cf. \cite{HT3}).
\par 
In \S \ref{section:preliminaries}, we briefly review the theory of framed curves and framed surfaces, see \cite{FT2, HT1} for details. 
We can treat space curves and surfaces with singular points by using the theory of framed curves and framed surfaces. 
Framed curves and framed surfaces are spaces curves and surfaces with moving frames, respectively. 
We introduce Bishop directions and Bishop frames along framed curves that play key roles in this paper. 
In \S \ref{section:parallel_curves_and_normal_surfaces}, we introduce parallel curves and normal surfaces of framed curves with respect to directions. 
Parallel curves of framed curves are direct generalizations of parallel curves of frontals in the Euclidean plane. 
Parallel curves of framed curves may not be framed curves. 
Actually, we give necessary and sufficient conditions in Proposition \ref{proposition:parallel_curves_condition}, see also \cite{HT3}. 
On the other hand, normal surfaces of framed curves are direct generalizations of principal normal surfaces of Frenet curves that are investigated in \cite{IT1}. 
Normal surfaces of framed curves may not be framed surfaces. 
Actually, we also give necessary and sufficient conditions in Proposition \ref{proposition:normal_surfaces_condition}. 
For normal surfaces of framed curves, we give characterizations of cross caps, cuspidal edges, swallowtails, and cuspidal cross caps, see Proposition \ref{proposition:cc} and Theorem \ref{proposition:ce_sw_ccr}. 
Furthermore, we introduce families of functions on framed curves that are useful for the study of normal surfaces and their striction curves. 
In \S \ref{sec_evo_invo}, we define circular evolutes and involutes of framed curves under suitable assumptions. 
We investigate properties of circular evolutes and involutes of framed curves. 
We can observe that taking circular evolutes and involutes of framed curves are opposite operations under suitable assumptions, similarly to evolutes and involutes of fronts in the Euclidean plane. 
Moreover, we give relations among singular points of normal surfaces, circular evolutes, and involutes of framed curves. 
We give examples to understand phenomena for normal surfaces, circular evolutes, and involutes of framed curves accordingly. 
\par 
All maps and manifolds considered here are differential of class $C^\infty$. 
\bigskip\\
\noindent
\textbf{Acknowledgement}. 
The second author was partially supported by JSPS KAKENHI Grant Number JP 20K03573. 
%
%
\section{Preliminaries}\label{section:preliminaries}
Let $\R^3$ be the 3-dimensional Euclidean space equipped with the canonical inner product $\bm{x} \cdot \bm{y} = x_1 y_1 + x_2 y_2 + x_3 y_3$, where $\bm{x} = (x_1, x_2, x_3)$ and $\bm{y} = (y_1, y_2, y_3)$. 
We define the norm of $\bm{x}$ by $\Vert \bm{x} \Vert = \sqrt{\bm{x} \cdot \bm{x}}$ and the vector product of $\bm{x}$ and $\bm{y}$ by 
\begin{equation*}
\bm{x} \times \bm{y} = \det 
\left( 
\begin{array}{ccc}
\bm{e}_1 & \bm{e}_2 & \bm{e}_3 \\
x_1 & x_2 & x_3 \\
y_1 & y_2 & y_3
\end{array}
\right), 
\end{equation*}
where $\{ \bm{e}_1, \bm{e}_2, \bm{e}_3 \}$ is the canonical basis of $\R^3$. 
We define the unit sphere in $\R^3$ by $S^2 = \{ \bm{x} \in \R^3 \ | \ \Vert \bm{x} \Vert =1 \}$ and the set of orthonormal vectors in $\R^3$ by $\Delta = \{ (\bm{x}, \bm{y}) \in S^2 \times S^2 \ | \ \bm{x} \cdot \bm{y} = 0 \}$. 
\par 
We briefly review the theory of framed curves. 
Framed curves are space curves with moving frames. 
We can treat space curves with degenerate points (namely, $\dot{\bm{\gamma}}(t) \times \ddot{\bm{\gamma}}(t) = 0$) and singular points (namely, $\dot{\bm{\gamma}}(t) = 0 $) by using the theory of framed curves, see \cite{HT1} for details. 
Here, $\dot{\bm{\gamma}}(t) = (d \bm{\gamma}/dt)(t)$ and $\ddot{\bm{\gamma}}(t) = (d^2 \bm{\gamma}/dt^2)(t)$. 
\begin{definition}
\rm 
We say that $(\bm{\gamma}, \bm{\nu}_1, \bm{\nu}_2):I \to \R^3 \times \Delta$ is a {\it framed curve} if $\dot{\bm{\gamma}}(t) \cdot \bm{\nu}_1(t) = 0$ and $\dot{\bm{\gamma}}(t) \cdot \bm{\nu}_2(t) = 0$ for all $t \in I$. 
Moreover, $\bm{\gamma}:I \to \R^3$ is a {\it framed base curve} if there exists $(\bm{\nu}_1, \bm{\nu}_2):I \to \Delta$ such that $(\bm{\gamma}, \bm{\nu}_1, \bm{\nu}_2)$ is a framed curve. 
\end{definition}
\par
Let $(\bm{\gamma}, \bm{\nu}_1, \bm{\nu}_2):I \to \R^3 \times \Delta$ be a framed curve. 
We define a (generalized) {\it unit tangent vector} $\bm{\mu}:I \to S^2$ of $\bm{\gamma}$ by $\bm{\mu}(t) = \bm{\nu}_1(t) \times \bm{\nu}_2(t)$. 
Then $\{ \bm{\nu}_1, \bm{\nu}_2, \bm{\mu} \}$ is a moving frame along $\bm{\gamma}$, and we have the following Frenet-Serret type formula: 
\begin{equation*}
\left(
\begin{array}{c}
\dot{\bm{\nu}}_1(t) \\
\dot{\bm{\nu}}_2(t) \\
\dot{\bm{\mu}}(t)
\end{array}
\right)
=
\left(
\begin{array}{ccc}
0 & \ell(t) & m(t) \\
-\ell(t) & 0 & n(t) \\
-m(t) & -n(t) & 0
\end{array}
\right)
\left(
\begin{array}{c}
\bm{\nu}_1(t) \\
\bm{\nu}_2(t) \\
\bm{\mu}(t)
\end{array}
\right), \quad \dot{\bm{\gamma}}(t) = \alpha(t) \bm{\mu}(t), 
\end{equation*}
where $\ell(t) = \dot{\bm{\nu}}_1(t) \cdot \bm{\nu}_2(t)$, $m(t) = \dot{\bm{\nu}}_1(t) \cdot \bm{\mu}(t)$, $n(t) = \dot{\bm{\nu}}_2(t) \cdot \bm{\mu}(t)$, and $\alpha(t) = \dot{\bm{\gamma}}(t) \cdot \bm{\mu}(t)$. 
We call the quadruple $(\ell, m, n, \alpha)$ a {\it curvature} of $(\bm{\gamma}, \bm{\nu}_1, \bm{\nu}_2)$. 
Note that $t_0$ is a singular point of $\bm{\gamma}$ if and only if $\alpha(t_0) = 0$. 
\begin{definition}
\rm 
Let $(\bm{\gamma}, \bm{\nu}_1, \bm{\nu}_2)$ and $(\widetilde{\bm{\gamma}}, \widetilde{\bm{\nu}}_1, \widetilde{\bm{\nu}}_2):I \to \R^3 \times \Delta$ be framed curves. 
We say that $(\bm{\gamma}, \bm{\nu}_1, \bm{\nu}_2)$ and $(\widetilde{\bm{\gamma}}, \widetilde{\bm{\nu}}_1, \widetilde{\bm{\nu}}_2)$ are {\it congruent} as framed curves if there exist a rotation $A \in SO(3)$ and a translation $\bm{a} \in \R^3$ such that $\widetilde{\bm{\gamma}}(t) = A(\bm{\gamma}(t)) + \bm{a}$ and $\widetilde{\bm{\nu}}_i(t) = A(\bm{\nu}_i(t)) \ (i = 1,2)$ for all $t \in I$. 
\end{definition}
\par 
We proved the following existence and uniqueness theorems for framed curves in \cite{HT1}, see also \cite{FT1}. 
\begin{theorem}[Existence theorem for framed curves, \cite{HT1}]
Let $(\ell, m, n, \alpha):I \to \R^4$ be a smooth mapping. 
Then there exists a framed curve $(\bm{\gamma}, \bm{\nu}_1, \bm{\nu}_2):I \to \R^3 \times \Delta$ whose curvature is given by $(\ell, m, n, \alpha)$. 
\end{theorem}
\begin{theorem}[Uniqueness theorem for framed curves, \cite{HT1}]
Let $(\bm{\gamma}, \bm{\nu}_1, \bm{\nu}_2)$ and $(\widetilde{\bm{\gamma}}, \widetilde{\bm{\nu}}_1, \widetilde{\bm{\nu}}_2):I \to \R^3 \times \Delta$ be framed curves with curvatures $(\ell, m, n, \alpha)$ and $(\widetilde{\ell}, \widetilde{m}, \widetilde{n}, \widetilde{\alpha})$, respectively. 
Then $(\bm{\gamma}, \bm{\nu}_1, \bm{\nu}_2)$ and $(\widetilde{\bm{\gamma}}, \widetilde{\bm{\nu}}_1, \widetilde{\bm{\nu}}_2)$ are congruent as framed curves if and only if the curvatures $(\ell, m, n, \alpha)$ and $(\widetilde{\ell}, \widetilde{m}, \widetilde{n}, \widetilde{\alpha})$ coincide. 
\end{theorem}
\par
For a framed curve $(\bm{\gamma}, \bm{\nu}_1, \bm{\nu}_2):I \to \R^3 \times \Delta$, we call a unit vector $\bm{v}(t) \in \langle \bm{\nu}_1(t), \bm{\nu}_2(t) \rangle_{\R}$ is a {\it Bishop direction} (or, {\it Bishop vector}) if there exists a smooth function $\beta:I \to \R$ such that $\dot{\bm{v}}(t) = \beta(t) \bm{\mu}(t)$, namely, $\bm{v}$ is a spherical frontal curve (cf. \cite{T1}). 
We also call a moving frame $\{ \bm{v}, \bm{w}, \bm{\mu} \}$ a {\it Bishop frame} along $\gamma$ if both $\bm{v}$ and $\bm{w}$ are Bishop directions (cf. \cite{Bishop}). 
Note that the original moving frame $\{ \bm{\nu}_1, \bm{\nu}_2, \bm{\mu} \}$ is a Bishop frame if and only if $\ell(t) = 0$ for all $t \in I$. 
\par 
We define $(\bm{v}, \bm{w}):I \to \Delta$ by 
\begin{equation*}
\left(
\begin{array}{c}
\bm{v}(t) \\
\bm{w}(t)
\end{array}
\right)
=
\left(
\begin{array}{cc}
\cos \theta(t) & - \sin \theta(t) \\
\sin \theta(t) & \cos \theta(t)
\end{array}
\right)
\left(
\begin{array}{c}
\bm{\nu}_1(t) \\
\bm{\nu}_2(t)
\end{array}	
\right), 
\end{equation*}
where $\theta:I \to \R$ is a smooth function. 
Then $(\bm{\gamma}, \bm{v}, \bm{w}):I \to \R^3 \times \Delta$ is also a framed curve with the unit tangent vector $\bm{\mu}$. 
The moving frame $\{ \bm{v}, \bm{w}, \bm{\mu} \}$ is called a {\it rotated frame} of $\{ \bm{\nu}_1, \bm{\nu}_2, \bm{\mu} \}$ by $\theta$. 
By straightforward calculations, we have 
\begin{align*}
\dot{\bm{v}}(t) &= \{\ell(t) - \dot{\theta}(t)\} \sin \theta(t) \bm{\nu}_1(t) + \{\ell(t) - \dot{\theta}(t)\} \cos \theta(t) \bm{\nu}_2(t) \\
& \qquad + \{m(t) \cos \theta(t) - n(t) \sin \theta(t)\} \bm{\mu}(t), \\
\dot{\bm{w}}(t) &= - \{ \ell(t) - \dot{\theta}(t) \} \cos \theta(t) \bm{\nu}_1(t) + \{ \ell(t) - \dot{\theta}(t) \} \sin \theta(t) \bm{\nu}_2(t) \\
& \qquad + \{ m(t) \sin \theta(t) + n(t) \cos \theta(t) \} \bm{\mu}(t). 
\end{align*}
Then the Frenet-Serret type formula of the rotated frame $\{ \bm{v}, \bm{w}, \bm{\mu} \}$ is given by 
\begin{equation}\label{curvature-rotated}
\left(
\begin{array}{c}
\dot{\bm{v}}(t) \\
\dot{\bm{w}}(t) \\
\dot{\bm{\mu}}(t)
\end{array}
\right)
=
\left(
\begin{array}{ccc}
0 & \overline{\ell}(t) & \overline{m}(t) \\
-\overline{\ell}(t) & 0 & \overline{n}(t) \\
- \overline{m}(t) & - \overline{n}(t) & 0
\end{array}
\right)
\left( 
\begin{array}{c}
\bm{v}(t) \\
\bm{w}(t) \\
\bm{\mu}(t)
\end{array}
\right), \quad \dot{\bm{\gamma}}(t) = \alpha(t) \bm{\mu}(t), 
\end{equation}
where $\overline{\ell}(t) = \ell(t) - \dot{\theta}(t)$, $\overline{m}(t) = m(t) \cos \theta(t) - n(t) \sin \theta(t)$, and $\overline{n}(t) = m(t) \sin \theta(t) + n(t) \cos \theta(t)$. 
If we take a function $\theta$ which satisfies $\dot{\theta}(t) = \ell(t)$ for all $t \in I$, then the rotated frame $\{ \bm{v}, \bm{w}, \bm{\mu} \}$ is a Bishop frame, so that we can construct a Bishop frame for any framed curve. 
\begin{remark}
\rm 
A Bishop frame $\{ \bm{v}, \bm{w}, \bm{\mu} \}$ isn't unique for a given framed curve $(\bm{\gamma}, \bm{\nu}_1, \bm{\nu}_2)$. 
We need only condition $\dot{\theta}(t) = \ell(t)$ for all $t \in I$, so that there is flexibility in integral constants. 
\end{remark}
\begin{example}\label{example:spherical_nephroid}
\rm 
Let $(\bm{\gamma}, \bm{\nu}_1, \bm{\nu}_2):[0, 2\pi) \to S^2 \times \Delta$ be 
\begin{align*}
\bm{\gamma}(t) &= \left( \frac{3}{4} \cos t - \frac{1}{4} \cos 3t, \frac{3}{4} \sin t - \frac{1}{4} \sin 3t, \frac{\sqrt{3}}{2} \cos t \right), \\
\bm{\nu}_1(t) &= \left( - \frac{3}{4} \sin t - \frac{1}{4} \sin 3t, \cos^3 t, \frac{\sqrt{3}}{2} \sin t \right), \\
\bm{\nu}_2(t) &= \left( \frac{3}{4} \cos t - \frac{1}{4} \cos 3t, \sin^3 t, \frac{\sqrt{3}}{2} \cos t \right). 
\end{align*}
Then $\bm{\gamma}$ is called the {\it spherical nephroid}, see Figure \ref{fig:spherical_nephroid}. 
By straightforward calculations, $(\bm{\gamma}, \bm{\nu}_1, \bm{\nu}_2)$ is a framed curve with the curvature $\ell(t) = 0$, $m(t) = - \sqrt{3} \cos t$, $n(t) = \sqrt{3} \sin t$, and $\alpha(t) = \sqrt{3} \sin t$, so that the original moving frame $\{ \bm{\nu}_1, \bm{\nu}_2, \bm{\mu} \}$ is a Bishop frame. 
\begin{figure}[htbp]
\begin{center}
\includegraphics[width=50mm]{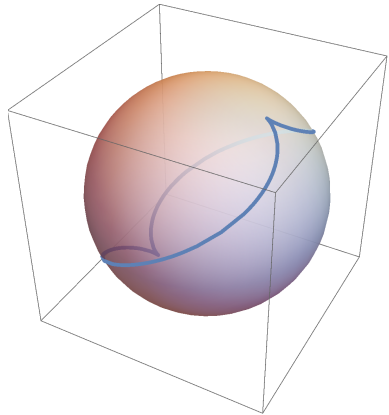}
\end{center}
\caption{The spherical nephroid.}
\label{fig:spherical_nephroid}
\end{figure}
\end{example}
\par
We briefly review the theory of framed surfaces. 
Framed surfaces are surfaces with moving frames. 
We can treat surfaces with singular points by using the theory of framed surfaces, see \cite{FT2} for details. 
Let $U$ be a simply connected domain of $\R^2$. 
\begin{definition}
\rm 
We say that $(\bm{x}, \bm{n}, \bm{s}):U \to \R^3 \times \Delta$ is a {\it framed surface} if $\bm{x}_u(u, v) \cdot \bm{n}(u, v) = 0$ and $\bm{x}_v(u, v) \cdot \bm{n}(u, v) = 0$ for all $(u, v) \in U$, where $\bm{x}_u(u, v) = (\partial \bm{x}/\partial u)(u, v)$ and $\bm{x}_v(u, v) = (\partial \bm{x}/\partial v)(u, v)$. 
Moreover, $\bm{x}:U \to \R^3$ is a {\it framed base surface} if there exists $(\bm{n}, \bm{s}):U \to \Delta$ such that $(\bm{x}, \bm{n}, \bm{s})$ is a framed surface. 
\end{definition}
\begin{definition}
\rm 
We say that $(\bm{x}, \bm{n}):U \to \R^3 \times S^2$ is a {\it Legendre surface} if $\bm{x}_u(u, v) \cdot \bm{n}(u, v) = 0$ and $\bm{x}_v(u, v) \cdot \bm{n}(u, v) = 0$ for all $(u, v) \in U$. 
We say that $(\bm{x}, \bm{n}):U \to \R^3 \times S^2$ is a {\it Legendre immersion} if $(\bm{x},\bm{n})$ is a Legendre surface and an immersion.  
Moreover, $\bm{x}:U \to \R^3$ is a {\it frontal} (respectively, {\it front}) if there exists $\bm{n}:U \to S^2$ such that $(\bm{x}, \bm{n})$ is a Legendre surface (respectively, Legendre immesion). 
\end{definition}
\par
Let \( (\bm{x}, \bm{n}, \bm{s}):U \to \R^3 \times \Delta \) be a framed surface. 
We define a (generalized) {\it tangent vector} \( \bm{t}:U \to S^2 \) of $\bm{x}$ by \( \bm{t}(u, v) = \bm{n}(u, v) \times \bm{s}(u, v) \). 
Then \( \{ \bm{n}, \bm{s}, \bm{t} \} \) is a moving frame along \( \bm{x} \), and we have the following structure equations: 
\begin{equation*}
\left(
\begin{array}{c}
\bm{x}_u(u, v) \\
\bm{x}_v(u, v)
\end{array}
\right)
=
\left(
\begin{array}{cc}
a_1(u,v) & b_1(u,v) \\
a_2(u,v) & b_2(u,v)
\end{array}
\right)
\left(
\begin{array}{c}
\bm{s}(u,v) \\
\bm{t}(u,v)
\end{array}
\right), 
\end{equation*}
\begin{equation*}
\left(
\begin{array}{c}
\bm{n}_u(u,v) \\
\bm{s}_u(u,v) \\
\bm{t}_u(u,v)
\end{array}
\right)
=
\left(
\begin{array}{ccc}
0 & e_1(u,v) & f_1(u,v) \\
- e_1(u,v) & 0 & g_1(u,v) \\
- f_1(u,v) & - g_1(u,v) & 0
\end{array}
\right)
\left(
\begin{array}{c}
\bm{n}(u,v) \\
\bm{s}(u,v) \\
\bm{t}(u,v)
\end{array}
\right), 
\end{equation*}
and
\begin{equation*}
\left(
\begin{array}{c}
\bm{n}_v(u,v) \\
\bm{s}_v(u,v) \\
\bm{t}_v(u,v)
\end{array}
\right)
=
\left(
\begin{array}{ccc}
0 & e_2(u,v) & f_2(u,v) \\
- e_2(u,v) & 0 & g_2(u,v) \\
- f_2(u,v) & - g_2(u,v) & 0
\end{array}
\right)
\left(
\begin{array}{c}
\bm{n}(u,v) \\
\bm{s}(u,v) \\
\bm{t}(u,v)
\end{array}
\right), 
\end{equation*}
where $a_i(u, v)$, $b_i(u, v)$, $e_i(u, v)$, $f_i(u, v)$, and $g_i(u, v)$ $(i = 1, 2)$ are smooth functions. 
We call the functions $(a_i, b_i, e_i, f_i, g_i)$ {\it basic invariants} of $(\bm{x}, \bm{n}, \bm{s})$. 
We represent the above matrices by
\begin{align*}
\mathcal{G}(u, v) = 
\left(
\begin{array}{cc}
a_1(u,v) & b_1(u,v) \\
a_2(u,v) & b_2(u,v)
\end{array}
\right), \quad \mathcal{F}_1(u, v) = 
\left(
\begin{array}{ccc}
0 & e_1(u,v) & f_1(u,v) \\
- e_1(u,v) & 0 & g_1(u,v) \\
- f_1(u,v) & - g_1(u,v) & 0
\end{array}
\right),
\end{align*}
and
\[
\mathcal{F}_2(u, v) = 
\left(
\begin{array}{ccc}
0 & e_2(u,v) & f_2(u,v) \\
- e_2(u,v) & 0 & g_2(u,v) \\
- f_2(u,v) & - g_2(u,v) & 0
\end{array}
\right).
\]
Note that $(u_0, v_0)$ is a singular point of $\bm{x}$ if and only if $\det \mathcal{G}(u_0, v_0) = 0$. 
\par
Since the integrability conditions $\bm{x}_{uv}(u, v) = \bm{x}_{vu}(u, v)$ and 
\begin{equation*}
\mathcal{F}_{2,u}(u, v) - \mathcal{F}_{1,v}(u, v) = \mathcal{F}_1(u, v) \mathcal{F}_2(u, v) - \mathcal{F}_2(u, v) \mathcal{F}_1(u, v), 
\end{equation*}
the basic invariants should be satisfied the following conditions: 
\begin{equation}\label{equation:integrability_condition_1}
\left\{
\begin{array}{rcl}
a_{1,v}(u, v) - b_1(u, v) g_2(u, v) & = & a_{2,u}(u, v) - b_2(u, v) g_1 (u, v), \\
b_{1,v}(u, v) - a_2(u, v) g_1(u, v) & = & b_{2,u}(u, v) - a_1(u, v) g_2(u, v), \\
a_1(u, v) e_2 (u, v) + b_1(u, v) f_2(u, v) & = & a_2(u, v) e_1(u, v) + b_2(u, v) f_1(u, v), 
\end{array}
\right.
\end{equation}
and
\begin{equation}\label{equation:integrability_condition_2}
\left\{
\begin{array}{rcl}
e_{1,v}(u, v) - f_1(u, v) g_2(u, v) & = & e_{2,u}(u, v) - f_2(u, v) g_1 (u, v), \\
f_{1,v}(u, v) - e_2(u, v) g_1(u, v) & = & f_{2,u}(u, v) - e_1(u, v) g_2(u, v), \\
g_{1,v}(u, v) - e_1(u, v) f_2(u, v) & = & g_{2,u}(u, v) - e_2(u, v) f_1(u, v) 
\end{array}
\right.
\end{equation}
for all $(u, v) \in U$. 
\begin{definition}
\rm 
Let $(\bm{x}, \bm{n}, \bm{s})$ and $(\widetilde{\bm{x}}, \widetilde{\bm{n}}, \widetilde{\bm{s}}):U \to \R^3 \times \Delta$ be framed surfaces. 
We say that $(\bm{x}, \bm{n}, \bm{s})$ and $(\widetilde{\bm{x}}, \widetilde{\bm{n}}, \widetilde{\bm{s}})$ are {\it congruent} as framed surfaces if there exist a rotation $A \in SO(3)$ and a translation $\bm{a} \in \R^3$ such that $\widetilde{\bm{x}}(u, v) = A(\bm{x}(u, v)) + \bm{a}$, $\widetilde{\bm{n}}(u, v) = A(\bm{n}(u, v))$, and $\widetilde{\bm{s}}(u, v) = A(\bm{s}(u, v))$ for all $(u, v) \in U$. 
\end{definition}
\par 
We have the following existence and uniqueness theorems for framed surfaces in \cite{FT2}. 
\begin{theorem}[Existence theorem for framed surfaces, \cite{FT2}]
Let $a_i(u, v)$, $b_i(u, v)$, $e_i(u, v)$, $f_i(u, v)$, and $g_i(u, v)$ $(i = 1, 2)$ be smooth functions with integrability conditions {\rm (\ref{equation:integrability_condition_1})} and {\rm (\ref{equation:integrability_condition_2})}. 
Then there exists a framed surface $(\bm{x}, \bm{n}, \bm{s}):U \to \R^3 \times \Delta$ whose basic invariants are given by $(a_i, b_i, e_i, f_i, g_i)$. 
\end{theorem}
\begin{theorem}[Uniqueness theorem for framed surfaces, \cite{FT2}]
Let $(\bm{x}, \bm{n}, \bm{s})$ and $(\widetilde{\bm{x}}, \widetilde{\bm{n}}, \widetilde{\bm{s}}):U \to \R^3 \times \Delta$ be framed surfaces with basic invariants $(\mathcal{G}, \mathcal{F}_1, \mathcal{F}_2)$ and $(\widetilde{\mathcal{G}}, \widetilde{\mathcal{F}}_1, \widetilde{\mathcal{F}}_2)$, respectively. 
Then $(\bm{x}, \bm{n}, \bm{s})$ and $(\widetilde{\bm{x}}, \widetilde{\bm{n}}, \widetilde{\bm{s}})$ are congruent as framed surfaces if and only if the basic invariants $(\mathcal{G}, \mathcal{F}_1, \mathcal{F}_2)$ and $(\widetilde{\mathcal{G}}, \widetilde{\mathcal{F}}_1, \widetilde{\mathcal{F}}_2)$ coincide. 
\end{theorem}
\par
For a framed surface $(\bm{x}, \bm{n}, \bm{s}):U \to \R^3 \times \Delta$, we define $C^F = (J^F, K^F, H^F):U \to \R^3$ by 
\begin{equation*}
J^F(u, v) = \det
\left(
\begin{array}{cc}
a_1(u, v) & b_1(u, v) \\
a_2(u, v) & b_2(u, v)
\end{array}
\right), 
\quad K^F(u, v) = \det
\left(
\begin{array}{cc}
e_1(u, v) & f_1(u, v) \\
e_2(u, v) & f_2(u, v)
\end{array}
\right), 
\end{equation*}
and
\begin{equation*}
H^F(u, v) = - \frac{1}{2}
\left(
\det
\left(
\begin{array}{cc}
a_1(u, v) & f_1(u, v) \\
a_2(u, v) & f_2(u, v)
\end{array}
\right)
-
\det
\left(
\begin{array}{cc}
b_1(u, v) & e_1(u, v) \\
b_2(u, v) & e_2(u, v)
\end{array}
\right)
\right).
\end{equation*}
We call $C^F$ a {\it curvature} of the framed surface $(\bm{x}, \bm{n}, \bm{s})$. 
The curvature is useful to recognize that the framed base surface is a front or not. 
\begin{proposition}[\cite{FT2, M1}]\label{pro:front_condition}
Let $(\bm{x}, \bm{n}, \bm{s}):U \to \R^3 \times \Delta$ be a framed surface. 
Then we have the following. 
\begin{enumerate}
\item[\rm (1)] $\bm{x}$ is a front around $(u_0, v_0)$ if and only if $C^F(u_0, v_0) \neq 0$. 
\item[\rm (2)] Suppose that $\rank \mathcal{G}(u_0, v_0) = 1$. 
Then $\bm{x}$ is a front around $(u_0, v_0)$ if and only if $H^F(u_0, v_0) \neq 0$. 
\end{enumerate}
\end{proposition}
%
%
\section{Parallel curves and normal surfaces}\label{section:parallel_curves_and_normal_surfaces}
In this section, we introduce parallel curves and normal surfaces of framed curves with respect to directions. 
Let $(\bm{\gamma}, \bm{\nu}_1, \bm{\nu}_2):I \to \R^3 \times \Delta$ be a framed curve with a rotated frame $\{ \bm{v}, \bm{w}, \bm{\mu} \}$. 
We denote the curvature of the framed curve $(\bm{\gamma},\bm{v},\bm{w}):I \to \R^3 \times \Delta$ by $(\overline{\ell},\overline{m},\overline{n},\alpha)$, see  \eqref{curvature-rotated}. 
\subsection{Parallel curves}
We introduce parallel curves of framed curves with respect to directions, see also \cite{HT2, HT3}. 
Then we define a space curve $P_{\bm{\gamma}}[\bm{v}]:I \to \R^3$ by
\begin{equation*}
P_{\bm{\gamma}}[\bm{v}](t) = \bm{\gamma}(t) + \lambda \bm{v}(t)
\end{equation*}
for a fixed $\lambda \in \R \setminus \{ 0 \}$. 
We call the space curve $P_{\bm{\gamma}}[\bm{v}]$ a {\it parallel curve} of $\bm{\gamma}$ with respect to $\bm{v}$ and $\lambda$ (or, $\bm{v}$-{\it parallel curve}). 
By a straightforward calculation, we have 
\begin{equation*}
\dot{P}_{\bm{\gamma}}[\bm{v}](t) = \lambda \overline{\ell}(t) \bm{w}(t) + \left\{ \alpha(t) + \lambda \overline{m}(t) \right\} \bm{\mu}(t). 
\end{equation*}
Therefore, the director $\bm{v}$ is one of the normal vectors of $P_{\bm{\gamma}}[\bm{v}]$. 
Note that $t_0$ is a singular point of $P_{\bm{\gamma}}[\bm{v}]$ if and only if $\overline{\ell}(t_0) = 0$ and $\alpha(t_0) + \lambda \overline{m}(t_0) = 0$. 
In order to consider the locus of singular values of $P_{\bm{\gamma}}[\bm{v}]$ continuously, we need a condition $\overline{\ell}(t) = 0$ for all $t \in I$. 
\par
A parallel curve of a framed curve may not be a framed base curve. 
Actually, we have the following necessary and sufficient conditions, see also \cite{HT3}. 
\begin{proposition}\label{proposition:parallel_curves_condition}
There exists a smooth map $\bm{n}:I \to S^2$ such that $(P_{\bm{\gamma}}[\bm{v}], \bm{v}, \bm{n}):I \to \R^3 \times \Delta$ is a framed curve if and only if there exists a smooth function $\varphi:I \to \R$ such that
\[
\lambda \overline{\ell}(t) \cos \varphi(t) - \left\{ \alpha(t) + \lambda \overline{m}(t) \right\} \sin \varphi(t) = 0
\]
for all $t \in I$ and a fixed $\lambda \in \R \setminus \{0\}$. 
\end{proposition}
\demo
Suppose that there exists a smooth map $\bm{n}:I \to S^2$ such that $(P_{\bm{\gamma}}[\bm{v}], \bm{v}, \bm{n})$ is a framed curve. 
Since $\{ \bm{v}, \bm{w}, \bm{\mu} \}$ is an orthonormal basis and $\bm{n}(t) \in \langle \bm{w}(t), \bm{\mu}(t) \rangle_{\R}$, there exists a smooth function $\varphi:I \to \R$ such that $\bm{n}(t) = \cos \varphi(t) \bm{w}(t) - \sin \varphi(t) \bm{\mu}(t)$. 
Then we have
\begin{equation*}
\dot{P}_{\bm{\gamma}}[\bm{v}](t) \cdot \bm{n}(t) = \lambda \overline{\ell}(t) \cos \varphi(t) - \{ \alpha(t) + \lambda \overline{m}(t) \} \sin \varphi(t) = 0
\end{equation*}
for all $t \in I$. 
\par 
Conversely, suppose that there exists a smooth function $\varphi:I \to \R$ such that $\lambda \overline{\ell}(t) \cos \varphi(t) - \left\{ \alpha(t) + \lambda \overline{m}(t) \right\} \sin \varphi(t) = 0$ for all $t \in I$. 
Then we define $\bm{n}:I \to S^2$ by $\bm{n}(t) = \cos \varphi(t) \bm{w}(t) - \sin \varphi(t) \bm{\mu}(t)$. 
By a straightforward calculation, $(P_{\bm{\gamma}}[\bm{v}], \bm{v}, \bm{n})$ is a framed curve. 
\enD
\begin{remark}\label{remark:bishop_parallel}
\rm If we take a Bishop frame $\{\bm{v}, \bm{w}, \bm{\mu} \}$, then $\overline{\ell}(t) = 0$ for all $t \in I$ and $P_{\bm{\gamma}}[\bm{v}]$ is always a framed base curve. 
More precisely, $(P_{\bm{\gamma}}[\bm{v}], \bm{v}, \bm{w})$ is a framed curve by Proposition \ref{proposition:parallel_curves_condition}. 
\end{remark}
\par 
Suppose that there exists a smooth function $\varphi:I \to \R$ such that
$$
\lambda \overline{\ell}(t) \cos \varphi(t) - \left\{ \alpha(t) + \lambda \overline{m}(t) \right\} \sin \varphi(t) = 0
$$
for all $t \in I$. 
We define $\bm{n}:I \to S^2$ by $\bm{n}(t) = \cos \varphi(t) \bm{w}(t) - \sin \varphi(t) \bm{\mu}(t)$. 
Then $(P_{\bm{\gamma}}[\bm{v}], \bm{v}, \bm{n})$ is a framed curve by Proposition \ref{proposition:parallel_curves_condition}. 
By straightforward calculations, we have the following curvature $(\ell_P, m_P, n_P, \alpha_P)$ of $(P_{\bm{\gamma}}[\bm{v}], \bm{v}, \bm{n})$: 
\begin{align*}
\ell_P(t) &= \dot{\bm{v}}(t) \cdot \bm{n}(t) = \overline{\ell}(t) \cos \varphi(t) - \overline{m}(t) \sin \varphi(t), \\
m_P(t) &= \dot{\bm{v}}(t) \cdot \bm{\mu}_P(t) = \overline{\ell}(t) \sin \varphi(t) + \overline{m}(t) \cos \varphi(t), \\
n_P(t) &= \dot{\bm{n}}(t) \cdot \bm{\mu}_P(t) = \overline{n}(t) - \dot{\varphi}(t),
\end{align*}
and
\begin{equation*}
\alpha_P(t) = \dot{P}_{\bm{\gamma}}[\bm{v}](t) \cdot \bm{\mu}_P(t) = \lambda \overline{\ell}(t) \sin \varphi(t) + \{ \alpha(t) + \lambda \overline{m}(t) \} \cos \varphi(t), 
\end{equation*}
where $\bm{\mu}_P(t) = \bm{v}(t) \times \bm{n}(t) = \sin \varphi(t) \bm{w}(t) + \cos \varphi(t) \bm{\mu}(t)$. 
%
%
\subsection{Normal surfaces}\label{sec_normal_surface}
We introduce normal surfaces of framed curves with respect to directions. 
Let $(\bm{\gamma}, \bm{\nu}_1, \bm{\nu}_2):I \to \R^3 \times \Delta$ be a framed curve with a rotated frame $\{ \bm{v}, \bm{w}, \bm{\mu} \}$. 
Then we define a surface $NS_{\bm{\gamma}}[\bm{v}]:I \times \R \to \R^3$ by
\begin{equation*}
NS_{\bm{\gamma}}[\bm{v}](t, \lambda) = \bm{\gamma}(t) + \lambda \bm{v}(t). 
\end{equation*}
We call the surface $NS_{\bm{\gamma}}[\bm{v}]$ a {\it normal surface} of $\bm{\gamma}$ with respect to $\bm{v}$ (or, {\it $\bm{v}$-normal surface}). 
The normal surface is a ruled surface, and we have 
\begin{equation*}
\det (\dot{\bm{\gamma}}(t), \bm{v}(t), \dot{\bm{v}}(t)) = \alpha(t) \overline{\ell}(t),
\end{equation*}
so that $NS_{\bm{\gamma}}[\bm{v}]$ is a developable surface if and only if $\alpha(t) \overline{\ell}(t) = 0$ on the regular part of $NS_{\bm{\gamma}}[\bm{v}]$ (cf. \cite{IT2}). 
By a straightforward calculation, we have
\begin{equation*}
NS_{\bm{\gamma}}[\bm{v}]_t (t, \lambda) \times NS_{\bm{\gamma}}[\bm{v}]_\lambda (t, \lambda) = \{ \alpha(t) + \lambda \overline{m}(t) \} \bm{w}(t) - \lambda \overline{\ell}(t) \bm{\mu}(t). 
\end{equation*}
Therefore, the director $\bm{v}$ is one of the tangent vectors of $NS_{\bm{\gamma}}[\bm{v}]$. 
Note that $(t_0, \lambda_0)$ is a singular point of $NS_{\bm{\gamma}}[\bm{v}]$ if and only if $\alpha(t_0) + \lambda_0 \overline{m}(t_0) = 0$ and $\lambda_0 \overline{\ell}(t_0) = 0$. 
Since $NS_{\bm{\gamma}}[\bm{v}]_\lambda(t, \lambda) = \bm{v}(t)$, $\rank ( d NS_{\bm{\gamma}}[\bm{v}])(t_0, \lambda_0) = 1$ at any singular point $(t_0, \lambda_0)$ of $NS_{\bm{\gamma}}[\bm{v}]$. 
\par 
A normal surface of a framed curve may not be a framed base surface. 
Actually, we have the following necessary and sufficient conditions. 
\begin{proposition}\label{proposition:normal_surfaces_condition}
There exists a smooth map $\bm{n}:I \times \R \to S^2$ such that $(NS_{\bm{\gamma}}[\bm{v}], \bm{n}, \bm{v}):I \times \R \to \R^3 \times \Delta$ is a framed surface if and only if there exists a smooth function $\varphi:I \times \R \to \R$ such that 
\begin{equation}\label{framed_surface_condition}
\lambda \overline{\ell}(t) \cos \varphi(t, \lambda) - \{ \alpha(t) + \lambda \overline{m}(t) \} \sin \varphi(t, \lambda) = 0
\end{equation}
for all $(t, \lambda) \in I \times \R$. 
\end{proposition}
\demo
Suppose that there exists a smooth map $\bm{n}:I \times \R \to S^2$ such that $(NS_{\bm{\gamma}}[\bm{v}], \bm{n}, \bm{v})$ is a framed surface. 
Since $\{ \bm{v}, \bm{w}, \bm{\mu} \}$ is an orthonormal basis and $\bm{n}(t, \lambda) \in \langle \bm{w}(t), \bm{\mu}(t) \rangle_{\R}$, there exists a smooth function $\varphi:I \times \R \to \R$ such that $\bm{n}(t, \lambda) = \cos \varphi(t, \lambda) \bm{w}(t) - \sin \varphi(t, \lambda) \bm{\mu}(t)$. 
Then we have
\begin{equation*}
NS_{\bm{\gamma}}[\bm{v}]_t (t, \lambda) \cdot \bm{n}(t, \lambda) = \lambda \overline{\ell}(t) \cos \varphi(t, \lambda) - \{ \alpha(t) + \lambda \overline{m}(t) \} \sin \varphi(t, \lambda) = 0
\end{equation*}
for all $(t, \lambda) \in I \times \R$. 
\par
Conversely, suppose that there exists a smooth function $\varphi:I \times \R \to \R$ such that 
\[
\lambda \overline{\ell}(t) \cos \varphi(t, \lambda) - \{ \alpha(t) + \lambda \overline{m}(t) \} \sin \varphi(t, \lambda) = 0
\]
for all $(t, \lambda) \in I \times \R$. 
Then we define $\bm{n}:I \times \R \to S^2$ by $\bm{n}(t, \lambda) = \cos \varphi(t, \lambda) \bm{w}(t) - \sin \varphi(t, \lambda) \bm{\mu}(t)$. 
By a straightforward calculation, $(NS_{\bm{\gamma}}[\bm{v}], \bm{n}, \bm{v})$ is a framed surface. 
\enD
\begin{remark}\label{bishop_normal}
\rm If we take a Bishop frame $\{\bm{v}, \bm{w}, \bm{\mu} \}$, then $\overline{\ell}(t) = 0$ for all $t \in I$ and $NS_{\bm{\gamma}}[\bm{v}]$ is always a framed base surface. 
More precisely, $(NS_{\bm{\gamma}}[\bm{v}], \bm{w}, \bm{v})$ is a framed surface by Proposition \ref{proposition:normal_surfaces_condition}. 
\end{remark}
\par
\begin{figure}[htbp]
\begin{center}
\includegraphics[width=50mm]{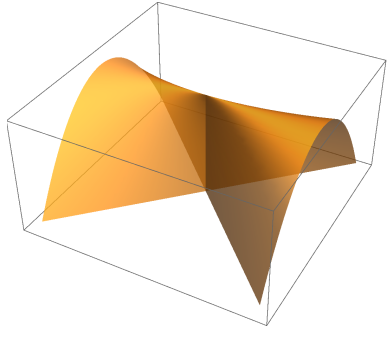}
\end{center}
\caption{A cross cap $(u, v) \mapsto (u, uv, v^2)$.}
\label{fig:cross_cap}
\end{figure}

 
We now consider cross caps of normal surfaces of framed curves. 
A singular point $p$ of a map $f:\R^2 \to \R^3$ is called a {\it cross cap} if the map-germ $f$ at $p$ is $\mathcal{A}$-equivalent to $(u, v) \mapsto (u, uv, v^2)$ at $0$ (cf. Figure \ref{fig:cross_cap}). 
If a singular point $p$ of $f$ is a cross cap, then $f$ is not a frontal at $p$. 
\par
By Proposition \ref{proposition:normal_surfaces_condition}, $NS_{\bm{\gamma}}[\bm{v}]$ may not be a framed base surface, so that $NS_{\bm{\gamma}}[\bm{v}]$ may have cross caps. 
\begin{proposition}\label{proposition:cc}
Suppose that $(t_0, \lambda_0) \in I \times \R$ is a singular point of $NS_{\bm{\gamma}}[\bm{v}]$, namely, $\lambda_0 \overline{\ell}(t_0) = 0$ and $\alpha(t_0) + \lambda_0 \overline{m}(t_0) = 0$. 
Then the singular point $(t_0, \lambda_0)$ of $NS_{\bm{\gamma}}[\bm{v}]$ is a cross cap if and only if 
$$
\alpha(t_0) \dot{\overline{\ell}}(t_0) + \dot{\alpha}(t_0) \overline{\ell}(t_0) \neq 0. 
$$
\end{proposition}
\demo
We apply Whitney's criterion for cross caps (cf. \cite{W1}). 
Since $NS_{\bm{\gamma}}[\bm{v}]_t(t_0, \lambda_0) = 0$ at a singular point $(t_0, \lambda_0)$, the singular point $(t_0, \lambda_0)$ of $NS_{\bm{\gamma}}[\bm{v}]$ is a cross cap if and only if $\det ( NS_{\bm{\gamma}}[\bm{v}]_\lambda, NS_{\bm{\gamma}}[\bm{v}]_{t\lambda}, NS_{\bm{\gamma}}[\bm{v}]_{tt})(t_0, \lambda_0) \neq 0$. 
By straightforward calculations, we have
\begin{align*}
NS_{\bm{\gamma}}[\bm{v}]_\lambda(t, \lambda) &= \bm{v}(t), \\
NS_{\bm{\gamma}}[\bm{v}]_{t\lambda}(t, \lambda) &= \overline{\ell}(t) \bm{w}(t) + \overline{m}(t) \bm{\mu}(t), \\
\end{align*}
and
\begin{align*}
NS_{\bm{\gamma}}[\bm{v}]_{tt}(t, \lambda) &= - \left\{ \alpha(t) \overline{m}(t) + \lambda \left( \overline{\ell}^2(t) + \overline{m}^2(t) \right) \right\} \bm{v}(t)  \\
& \quad - \left\{ \alpha(t) \overline{n}(t) + \lambda \left( \overline{m}(t) \overline{n}(t) - \dot{\overline{\ell}}(t) \right)  \right\} \bm{w}(t) \\
& \quad + \left\{ \dot{\alpha}(t) + \lambda \left( \overline{\ell}(t) \overline{n}(t) + \dot{\overline{m}}(t) \right) \right\} \bm{\mu}(t). 
\end{align*}
Therefore, 
\[
\det ( NS_{\bm{\gamma}}[\bm{v}]_\lambda, NS_{\bm{\gamma}}[\bm{v}]_{t\lambda}, NS_{\bm{\gamma}}[\bm{v}]_{tt})(t_0, \lambda_0) = \alpha(t_0) \dot{\overline{\ell}}(t_0) + \dot{\alpha}(t_0) \overline{\ell}(t_0). 
\]
This completes the proof. 
\enD
\par
\begin{remark}
\rm By Proposition \ref{proposition:cc}, if we take a Bishop vector $\bm{v}$ (namely, $\overline{\ell}(t) = 0$ for all $t \in I$), then singular points of $NS_{\bm{\gamma}}[\bm{v}]$ are never cross caps. 
\end{remark}
\begin{example}\label{example:astroid}
\rm Let $(\bm{\gamma}, \bm{\nu}_1, \bm{\nu}_2):[0, 2\pi) \to \R^3 \times \Delta$ be 
\begin{align*}
\bm{\gamma}(t) = \left( \cos^3 t, \sin^3 t, \cos 2t \right), \ \bm{\nu}_1(t) = \left( - \sin t, - \cos t, 0 \right), \ \bm{\nu}_2(t) = \left(-\frac{4}{5} \cos t, \frac{4}{5} \sin t, \frac{3}{5} \right). 
\end{align*}
Then $\bm{\gamma}$ is called the {\it spatial astroid}, see Figure \ref{fig:spatial_astroid}. 
By straightforward calculations, $(\bm{\gamma}, \bm{\nu}_1, \bm{\nu}_2)$ is a framed curve with the curvature $\ell(t) = 4/5$, $m(t) = 3/5$, $n(t) = 0,$ and $\alpha(t) = 5 \cos t \sin t$. 
If we take a rotated frame $\{ \bm{v}, \bm{w}, \bm{\mu} \}$ by $\theta(t) = 0$ (namely, $\bm{v}(t) = \bm{\nu}_1(t)$ and $\bm{w}(t) = \bm{\nu}_2(t)$), then $NS_{\bm{\gamma}}[\bm{v}]$ is given by 
$$
NS_{\bm{\gamma}}[\bm{v}](t,\lambda) = ( \cos^3 t - \lambda \sin t, \sin^3 t - \lambda \cos t, \cos 2t ). 
$$
By Proposition \ref{proposition:cc}, singular points $(k \pi/2,0)$ $(k = 0, 1, 2, 3)$ of $NS_{\bm{\gamma}}[\bm{v}]$ are cross caps, see Figure \ref{fig:spatial_astroid}. 
\begin{figure}[htbp]
\begin{minipage}{0.5\hsize}
\begin{center}
\includegraphics[width=50mm]{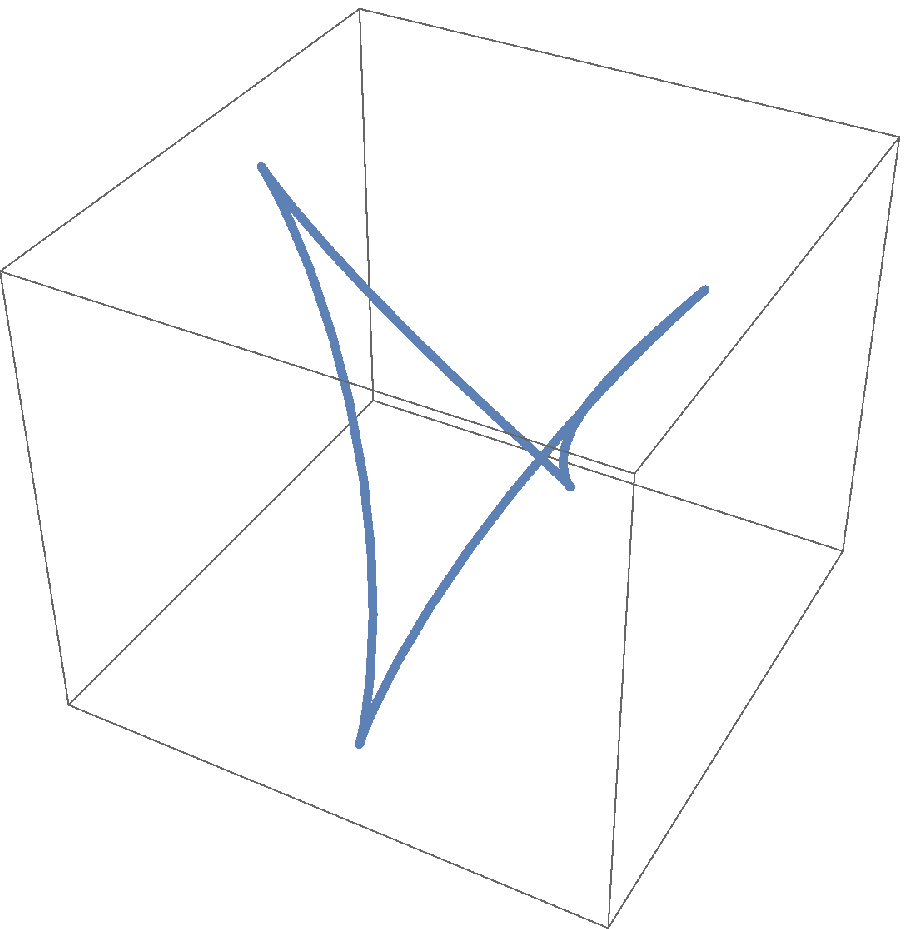}
\end{center}
\end{minipage}
\begin{minipage}{0.5\hsize}
\begin{center}
\includegraphics[width=50mm]{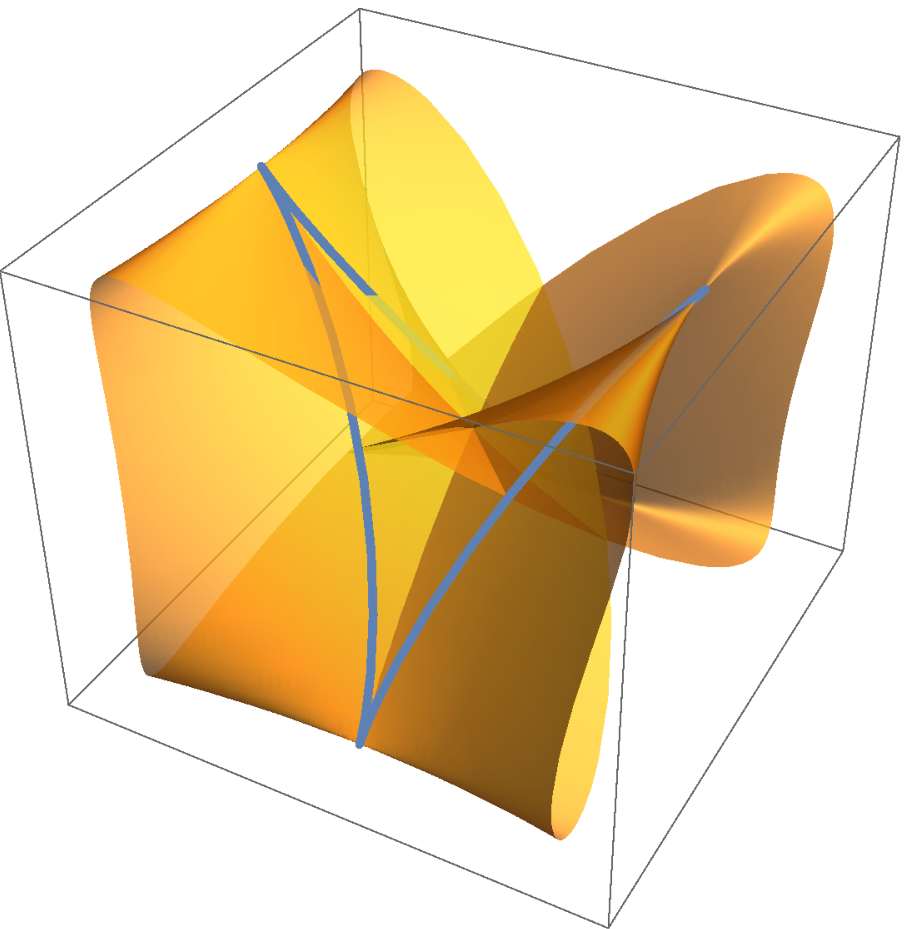}
\end{center}
\end{minipage}
\caption{Left to right: $\bm{\gamma}$ and $(\bm{\gamma}, NS_{\bm{\gamma}}[\bm{v}])$ of Example \ref{example:astroid}.}
\label{fig:spatial_astroid}
\end{figure}
\end{example}
\par 
Suppose that there exists a smooth function $\varphi:I \times \R \to \R$ such that 
\[
\lambda \overline{\ell}(t) \cos \varphi(t, \lambda) - \{ \alpha(t) + \lambda \overline{m}(t) \} \sin \varphi(t, \lambda) = 0
\]
for all $(t, \lambda) \in I \times \R$. 
We define $\bm{n}:I \times \R \to S^2$ by $\bm{n}(t,\lambda) = \cos \varphi(t, \lambda) \bm{w}(t) - \sin \varphi(t, \lambda) \bm{\mu}(t)$. 
Then $(NS_{\bm{\gamma}}[\bm{v}], \bm{n}, \bm{v})$ is a framed surface by Proposition \ref{proposition:normal_surfaces_condition}. 
By straightforward calculations, we have the following basic invariants $(a_i, b_i, e_i, f_i, g_i)$ $(i = 1, 2)$ of $(NS_{\bm{\gamma}}[\bm{v}], \bm{n}, \bm{v})$: 
\begin{align*}
a_1(t, \lambda) &= 0, \\
b_1(t, \lambda) &= - \lambda \overline{\ell}(t) \sin \varphi(t, \lambda) - \{ \alpha(t) + \lambda \overline{m}(t) \} \cos \varphi(t, \lambda), \\
e_1(t, \lambda) &= - \overline{\ell}(t) \cos \varphi(t, \lambda) + \overline{m}(t) \sin \varphi(t, \lambda), \\
f_1(t, \lambda) &= \varphi_t(t, \lambda) - \overline{n}(t), \\
g_1(t, \lambda) &= - \overline{\ell}(t) \sin \varphi(t, \lambda) - \overline{m}(t) \cos \varphi(t, \lambda), \\
a_2(t, \lambda) &= 1, \\
b_2(t, \lambda) &= 0, \\
e_2(t, \lambda) &= 0, \\
f_2(t, \lambda) &= \varphi_\lambda(t, \lambda), \\
g_2(t, \lambda) &= 0. 
\end{align*}
Naturally, the above basic invariants satisfy integrability conditions (\ref{equation:integrability_condition_1}) and (\ref{equation:integrability_condition_2}) under assumption (\ref{framed_surface_condition}), so that we have 
\begin{equation*}
\left\{
\begin{array}{rcl}
b_{1,\lambda}(t, \lambda) &=& g_1(t,\lambda), \\
b_1(t, \lambda) f_2(t, \lambda) &=& e_1(t, \lambda), \\
e_{1,\lambda}(u, v) &=& - f_2(u, v) g_1(u, v), \\
g_{1,\lambda}(u, v) &=& e_1(u, v) f_2(u, v). 
\end{array}
\right.
\end{equation*}
By the above integrability conditions, if $\lambda_0 \overline{\ell}(t_0) = 0$ and $\alpha(t_0) + \lambda_0 \overline{m}(t_0) = 0$, then $\alpha(t_0) \sin \varphi(t_0) = 0$. 
By the straightforward calculations, we have the following curvature $C^F = (J^F, K^F, H^F)$ of the framed surface $(NS_{\bm{\gamma}}[\bm{v}], \bm{n}, \bm{v})$: 
\begin{align*}
J^F(t, \lambda) &= \lambda \overline{\ell}(t) \sin \varphi(t, \lambda) + \{ \alpha(t) + \lambda \overline{m}(t) \} \cos \varphi(t, \lambda), \\
K^F(t, \lambda) &= \{ - \overline{\ell}(t) \cos \varphi(t, \lambda) + \overline{m}(t) \sin \varphi(t, \lambda) \} \varphi_\lambda(t, \lambda), \\
H^F(t, \lambda) &= \frac{1}{2} \{ \varphi_t(t, \lambda) - \overline{n}(t) \}. 
\end{align*}
By Proposition \ref{pro:front_condition} (2),  $NS_{\bm{\gamma}}[\bm{v}]$ is a front around a singular point $(t_0, \lambda_0) \in I \times \R$ if and only if $b_1(t_0, \lambda_0) = 0$ and $H^F(t_0, \lambda_0) \neq 0$. 
Under assumption (\ref{framed_surface_condition}), $b_1(t_0, \lambda_0) = 0$ if and only if $\lambda_0 \overline{\ell}(t_0) = 0$ and $\alpha(t_0) + \lambda_0 \overline{m}(t_0) = 0$. 
\par
\begin{figure}[htbp]
\begin{minipage}{0.325\hsize}
\begin{center}
\includegraphics[width=50mm]{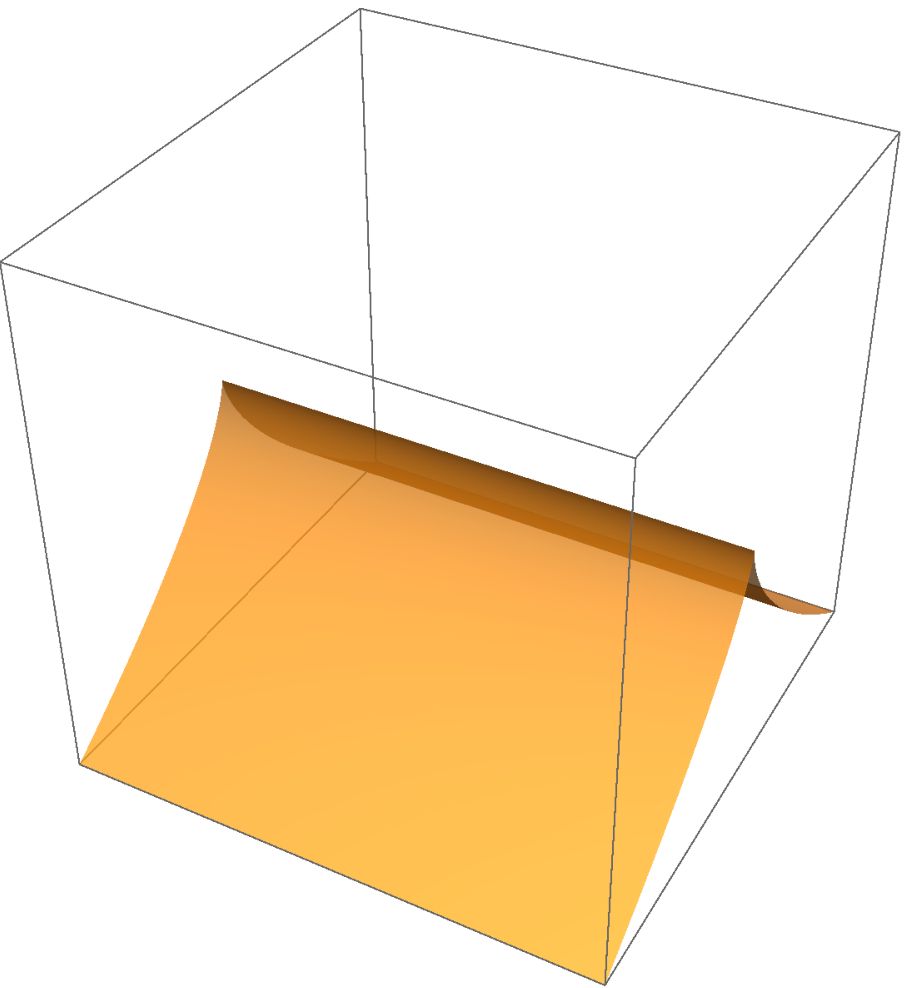}
\end{center}
\end{minipage}
\begin{minipage}{0.325\hsize}
\begin{center}
\includegraphics[width=50mm]{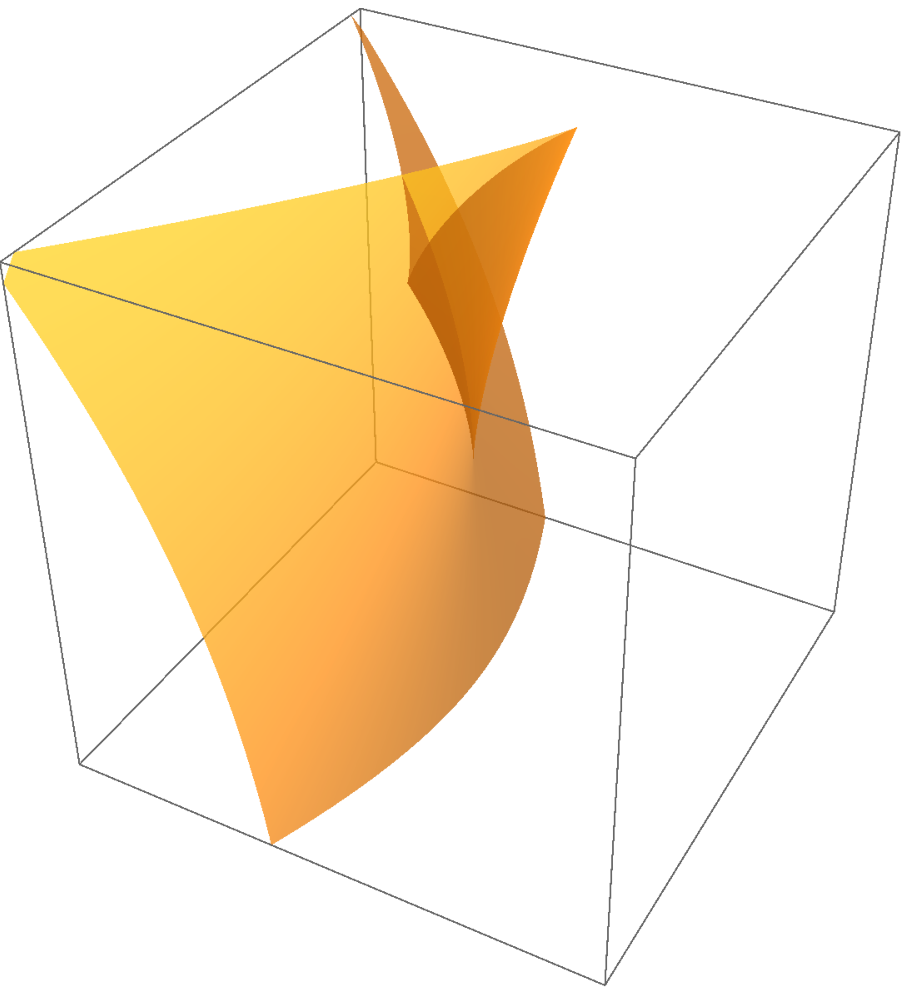}
\end{center}
\end{minipage}
\begin{minipage}{0.325\hsize}
\begin{center}
\includegraphics[width=50mm]{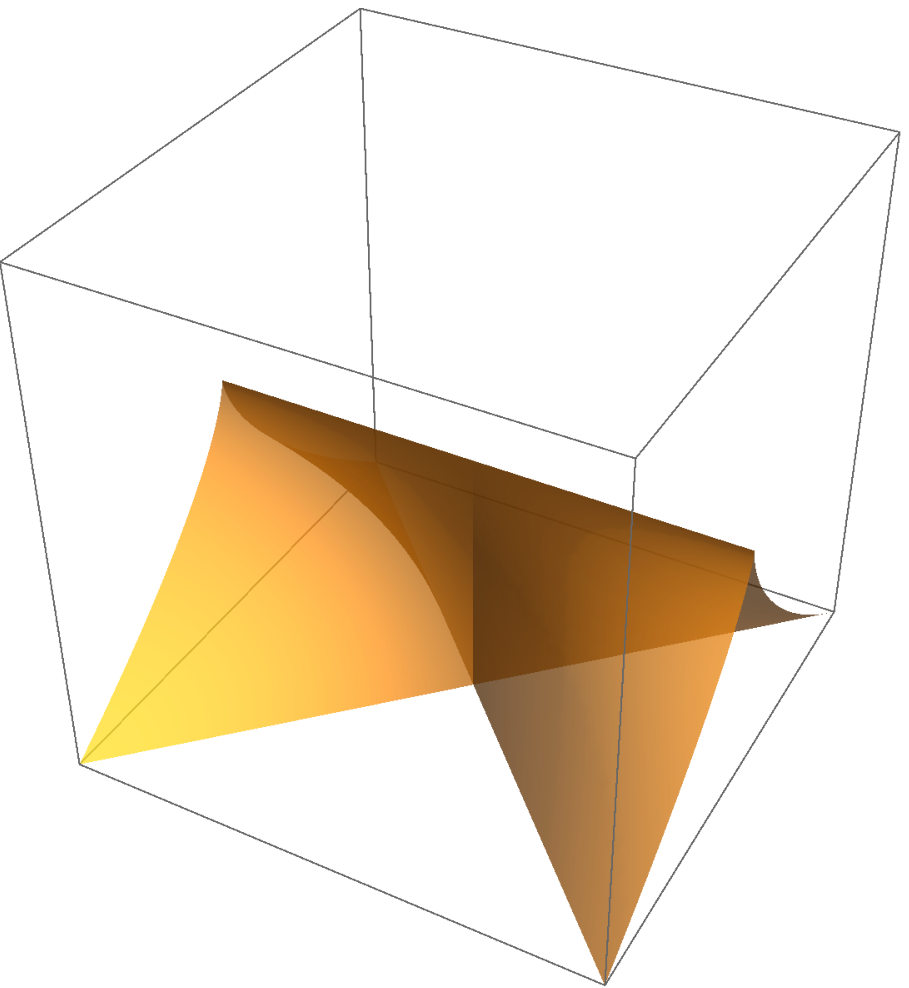}
\end{center}
\end{minipage}
\caption{Left to right: a cuspidal edge $(u, v) \mapsto (u, v^2, v^3)$, a swallowtail $(u, v) \mapsto (3u^4 + u^2 v, 4u^3 + 2uv, v)$,  and a cuspidal cross cap $(u,v) \mapsto (u, uv^3, v^2)$.}
\label{fig:singularities}
\end{figure}

We now consider cuspidal edges, swallowtails, and cuspidal cross caps of normal surfaces of framed curves. 
A singular point $p$ of a map $f:\R^2 \to \R^3$ is called a {\it cuspidal edge} (respectively, {\it swallowtail} and {\it cuspidal cross cap}) if the map-germ $f$ at $p$ is $\mathcal{A}$-equivalent to $(u, v) \mapsto (u, v^2, v^3)$ (respectively, $(u, v) \mapsto (3u^4 + u^2 v, 4u^3 + 2uv, v)$ and $(u,v) \mapsto (u, uv^3, v^2)$) at $0$, see Figure \ref{fig:singularities}. 
If a singular point $p$ of $f$ is a cuspidal edge or a swallowtail, then $f$ is a front at $p$. 
If a singular point $p$ of $f$ is a cuspidal cross cap, then $f$ is a frontal (not front) at $p$. 
\par 
We consider the signed area density function $\Lambda(t, \lambda) = \det (NS_{\bm{\gamma}}[\bm{v}]_t, NS_{\bm{\gamma}}[\bm{v}]_\lambda, \bm{n})(t, \lambda) = - b_1(t, \lambda)$.  
We denote the set of singular points of $NS_{\bm{\gamma}}[\bm{v}]$ by $S(NS_{\bm{\gamma}}[\bm{v}])$. 
Let $(t_0, \lambda_0) \in S(NS_{\bm{\gamma}}[\bm{v}])$ be a non-degenerate singular point of $NS_{\bm{\gamma}}[\bm{v}]$ (namely, $d \Lambda(t_0, \lambda_0) \neq 0$).  
Then there exists a regular curve $\bm{c}:(-\varepsilon, \varepsilon) \to I \times \R, s \mapsto \bm{c}(s)$ $(\varepsilon > 0)$ which gives a local parametrization of $S(NS_{\bm{\gamma}}[\bm{v}])$ by the implicit function theorem. 
Furthermore, there exists a null vector field $\bm{\eta}$ on $S(NS_{\bm{\gamma}}[\bm{v}])$ such that $dNS_{\bm{\gamma}}[\bm{v}](\bm{\eta}(s)) = 0$. 
We can take $\bm{\eta}$ as $\bm{\eta}(s) = (1, 0)$ for $NS_{\bm{\gamma}}[\bm{v}]$. 
We define a function $\Phi:(-\varepsilon, \varepsilon) \to \R$ by $\Phi(s) = \det ((NS_{\bm{\gamma}}[\bm{v}] \circ \bm{c})', \bm{v} \circ \bm{c}, d \bm{v}(\bm{\eta}))(s)$. 
By using the criteria in \cite{FSUY1, KRSUY1}, we give the following criteria of singularities of normal surfaces of framed curves. 
\begin{theorem}\label{proposition:ce_sw_ccr}
Let $(\bm{\gamma}, \bm{\nu}_1, \bm{\nu}_2):I \to \R^3 \times \Delta$ be a framed curve with a rotated frame $\{ \bm{v}, \bm{w}, \bm{\mu} \}$ and condition {\rm (\ref{framed_surface_condition})}. 
Suppose that $(t_0, \lambda_0) \in S(NS_{\bm{\gamma}}[\bm{v}])$ is a non-degenerate singular point of $NS_{\bm{\gamma}}[\bm{v}]$. 
\begin{itemize}
\item[\rm (1)] Suppose that $NS_{\bm{\gamma}}[\bm{v}]$ is a front at $(t_0, \lambda_0)$, namely, $H^F(t_0, \lambda_0) \neq 
0$.
Then we have the following. 
	\begin{enumerate}
	\item[\rm (i)] The singular point $(t_0, \lambda_0)$ of $NS_{\bm{\gamma}}[\bm{v}]$ is a cuspidal edge if and only if $b_{1,t}(t_0,\lambda_0) \neq 0$. 
	\item[\rm (ii)] The singular point $(t_0, \lambda_0)$ of $NS_{\bm{\gamma}}[\bm{v}]$ is a swallowtail if and only if $b_{1, t}(t_0,\lambda_0) = 0$ and $b_{1, tt}(t_0,\lambda_0) \neq 0$. 
	\end{enumerate}
\item[\rm (2)] Suppose that $NS_{\bm{\gamma}}[\bm{v}]$ is a frontal (not a front) at $(t_0, \lambda_0)$, namely, $C^F(t_0, \lambda_0) = 
0$. 
Then the singular point $(t_0, \lambda_0)$ of $NS_{\bm{\gamma}}[\bm{v}]$ is a cuspidal cross cap if and only if $b_{1,t}(t_0, \lambda_0) \neq 0$ and $b_{1, t}(t_0, \lambda_0) f_{1, \lambda}(t_0, \lambda_0) - b_{1, \lambda}(t_0, \lambda_0) f_{1, t}(t_0, \lambda_0) \neq 0$. 
\end{itemize}
Here, 
\begin{align*}
b_{1, t}(t_0, \lambda_0) =& - \left( \dot{\alpha}(t_0) + \lambda_0 \dot{\overline{m}}(t_0) \right) \cos \varphi(t_0, \lambda_0) - \lambda_0 \dot{\overline{\ell}}(t_0) \sin \varphi(t_0, \lambda_0), \\
b_{1, tt}(t_0, \lambda_0) =& - \left\{ \left( \ddot{\alpha}(t_0) + \lambda_0 \ddot{\overline{m}}(t_0) \right) + 2 \lambda_0 \dot{\overline{\ell}}(t_0) \varphi_t(t_0, \lambda_0) \right\} \cos \varphi(t_0, \lambda_0) \\
& - \left\{ \lambda_0 \ddot{\overline{\ell}}(t_0) - 2 \left( \dot{\alpha}(t_0)  + \lambda_0 \dot{\overline{m}}(t_0) \right) \varphi_t(t_0, \lambda_0) \right\} \sin \varphi(t_0, \lambda_0), 
\end{align*}
and
\begin{align*}
& b_{1, t}(t_0, \lambda_0) f_{1, \lambda}(t_0, \lambda_0) - b_{1, \lambda}(t_0, \lambda_0) f_{1, t}(t_0, \lambda_0) \\
& \quad = - \left\{ \left( \dot{\alpha}(t_0) + \lambda_0 \dot{\overline{m}}(t_0) \right) \varphi_{t \lambda}(t_0, \lambda_0) - \overline{m}(t_0) \left( \varphi_{t t}(t_0, \lambda_0) - \dot{\overline{n}}(t_0) \right) \right\} \cos \varphi(t_0, \lambda_0) \\
& \qquad - \left\{ \lambda_0 \overline{\ell}(t_0) \varphi_{t \lambda}(t_0, \lambda_0) - \overline{\ell}(t_0) \left( \varphi_{t t}(t_0, \lambda_0) - \dot{\overline{n}}(t_0) \right) \right\} \sin \varphi(t_0, \lambda_0). 
\end{align*}
\end{theorem}
\demo
(1) Suppose that $d \Lambda(t_0, \lambda_0) \neq 0$ and $H^F(t_0, \lambda_0) \neq 0$. 
By using criteria for cuspidal edges and swallowtails in \cite{KRSUY1}, we have assertions (i) and (ii) immediately. 
\par
(2) Suppose that $d \Lambda(t_0, \lambda_0) \neq 0$ and $C_F(t_0, \lambda_0) = 0$. 
Then we have $f_1(t_0, \lambda_0) = 0$. 
We divide the proof into two cases. 
\par 
(a) If $\Lambda_\lambda(t_0, \lambda_0) = b_{1, \lambda}(t_0, \lambda_0) \neq 0$, then there exists a smooth function $\psi:(I, t_0) \to (\R, \lambda_0)$ such that $\Lambda(t, \psi(t)) = 0$ and 
\[
\frac{d \psi}{dt}(t) = - \left. \frac{\Lambda_t(t, \lambda)}{\Lambda_\lambda(t, \lambda)} \right|_{\lambda = \psi(t)}
\]
by the implicit function theorem. 
Moreover, we have a regular curve $\bm{c}(t) = (t, \psi(t))$ and a null vector field $\bm{\eta}(t) = (1, 0)$. 
Then $\bm{\eta}\Lambda(t_0, \lambda_0) \neq 0$ is equivalent to $b_{1, t}(t_0, \lambda_0) \neq 0$. 
By straightforward calculations, we have
\begin{align*}
\left( NS_{\bm{\gamma}}[\bm{v}] \circ \bm{c} \right)' (t) &= \dot{\psi}(t) \bm{v}(t) + \psi(t) \overline{\ell}(t) \bm{w}(t)  + \left( \alpha(t) + \psi(t) \overline{m}(t) \right) \bm{\mu}(t), \\
\left( \bm{n} \circ \bm{c} \right)(t) &= \cos \varphi(t, \psi(t)) \bm{w}(t) - \sin \varphi(t, \psi(t)) \bm{\mu}(t), \\
d \bm{n}(\bm{\eta})(t) &= e_1(t, \psi(t)) \bm{v}(t) - f_1(t, \psi(t)) \sin \varphi(t, \psi(t)) \bm{w}(t) - f_1(t, \psi(t)) \cos \varphi(t, \psi(t)) \bm{\mu}(t), 
\end{align*}
and
\begin{align*}
\Phi(t) &= \det \left( \left( NS_{\bm{\gamma}}[\bm{v}] \circ \bm{c} \right)' (t), \left( \bm{n} \circ \bm{c} \right)(t), d \bm{n}(\bm{\eta})(t) \right) \\
&= - \dot{\psi}(t) f_1(t, \psi(t)) + b_1(t, \psi(t)) e_1(t, \psi(t)). 
\end{align*}
It follows that $\Phi(t_0) = 0$ and 
\[
\dot{\Phi}(t_0) = -\frac{b_{1,t}(t_0, \lambda_0) \{ b_{1, t}(t_0, \lambda_0) f_{1, \lambda}(t_0, \lambda_0) - b_{1, \lambda}(t_0, \lambda_0) f_{1, t}(t_0, \lambda_0) \}}{(b_{1, \lambda}(t_0, \lambda_0))^2}. 
\]
By using the criterion for cuspidal cross caps in \cite{FSUY1}, the singular point $(t_0, \lambda_0)$ of $NS_{\bm{\gamma}}[\bm{v}]$ is a cuspidal cross cap if and only if $b_{1,t}(t_0, \lambda_0) \neq 0$, $b_{1, \lambda}(t_0, \lambda_0) \neq 0$, and $b_{1, t}(t_0, \lambda_0) f_{1, \lambda}(t_0, \lambda_0) - b_{1, \lambda}(t_0, \lambda_0) f_{1, t}(t_0, \lambda_0) \neq 0$. 
\par 
(b) If $\Lambda_t(t_0, \lambda_0) = b_{1, t}(t_0, \lambda_0) \neq 0$ and $\Lambda_\lambda(t_0, \lambda_0) = b_{1, \lambda}(t_0, \lambda_0) = 0$, then there exists a smooth function $\overline{\psi}:(\R, \lambda_0) \to (\R, t_0)$ such that $\Lambda(\overline{\psi}(\lambda), \lambda) = 0$ and
\[
\frac{d \overline{\psi}}{d \lambda}(\lambda) = - \left. \frac{\Lambda_\lambda(t, \lambda)}{\Lambda_t(t, \lambda)} \right|_{t = \overline{\psi}(\lambda)}
\]
by the implicit function theorem. 
Moreover, we have a regular curve $\overline{\bm{c}}(\lambda) = ( \overline{\psi}(\lambda), \lambda )$ and a null vector field $\bm{\eta}(\lambda)  = (1, 0)$. 
Then $\bm{\eta} \Lambda(t_0, \lambda_0) \neq 0$ is equivalent to $b_{1, t} (t_0, \lambda_0) \neq 0$. 
By straightforward calculations, we have
\begin{align*}
\left( NS_{\bm{\gamma}}[\bm{v}] \circ \overline{\bm{c}} \right)' (\lambda) &= \bm{v}(\overline{\psi}(\lambda)) + \lambda \dot{\overline{\psi}}(\lambda) \overline{\ell}(\overline{\psi}(\lambda)) \bm{w}(\overline{\psi}(\lambda)) + \left( \alpha(\overline{\psi}(\lambda)) + \lambda \overline{m}(\overline{\psi}(\lambda)) \right) \dot{\overline{\psi}}(\lambda) \bm{\mu}(\overline{\psi}(\lambda)), \\
\left( \bm{n} \circ \overline{\bm{c}} \right) (\lambda) &= \cos \varphi(\overline{\psi}(\lambda), \lambda) \bm{w}(\overline{\psi}(\lambda)) - \sin \varphi(\overline{\psi}(\lambda), \lambda) \bm{\mu}(\overline{\psi}(\lambda)),\\
d \bm{n}(\bm{\eta})(\lambda) &= e_1(\overline{\psi}(\lambda), \lambda) \bm{v}(\overline{\psi}(\lambda)) - f_1(\overline{\psi}(\lambda), \lambda) \sin \varphi(\overline{\psi}(\lambda), \lambda) \bm{w}(\overline{\psi}(\lambda)) \\
& \qquad - f_1(\overline{\psi}(\lambda), \lambda) \cos \varphi(\overline{\psi}(\lambda), \lambda) \bm{\mu}(\overline{\psi}(\lambda)), 
\end{align*}
and
\begin{align*}
\overline{\Phi}(\lambda) &= \det \left( \left( NS_{\bm{\gamma}}[\bm{v}] \circ \overline{\bm{c}} \right)' (\lambda), \left( \bm{n} \circ \overline{\bm{c}} \right) (\lambda), d \bm{n}(\bm{\eta})(\lambda) \right) \\
&= - f_1(\overline{\psi}(\lambda), \lambda) + \dot{\overline{\psi}}(\lambda) e_1(\overline{\psi}(\lambda), \lambda) b_1(\overline{\psi}(\lambda),\lambda). 
\end{align*}
It follows that $\overline{\Phi}(\lambda_0) = 0$ and $\dot{\overline{\Phi}}(\lambda_0) = - f_{1, \lambda}(t_0, \lambda_0)$. 
By using the criterion for cuspidal cross caps in \cite{FSUY1}, the singular point $(t_0, \lambda_0)$ of $NS_{\bm{\gamma}}[\bm{v}]$ is a cuspidal cross cap if and only if $b_{1, t} (t_0, \lambda_0) \neq 0$ and $f_{1, \lambda}(t_0, \lambda_0) \neq 0$. 
\par
By summarizing cases (a) and (b), we have assertion (2). 
This completes the proof. 
\enD
\subsection{Normal surfaces of Bishop directions}\label{section:normal_developable_surfaces}
We now consider normal surfaces of framed curves with respect to Bishop directions. 
Let $(\bm{\gamma}, \bm{\nu}_1, \bm{\nu}_2):I \to \R^3 \times \Delta$ be a framed curve with a Bishop frame $\{ \bm{v}, \bm{w}, \bm{\mu} \}$. 
Then we have $\overline{\ell}(t) = 0$ and $\det \left( \dot{\bm{\gamma}}(t), \bm{v}(t), \dot{\bm{v}}(t) \right) = 0$ for all $t \in I$, so that $NS_{\bm{\gamma}}[\bm{v}]$ is a developable surface. 
We remark that developable surfaces are classified into cylinders, cones, and tangent surfaces (cf. \cite{V1}). 
$NS_{\bm{\gamma}}[\bm{v}]$ is a cylinder (namely, the director $\bm{v}$ is constant) if and only if $\overline{m}(t) = 0$ for all $t \in I$. 
Suppose that $NS_{\bm{\gamma}}[\bm{v}]$ is non-cylindrical, namely, $\overline{m}(t) \neq 0$ for all $t \in I$. 
Then the striction curve $\bm{\sigma}:I \to \R^3$ is given by
$$
\bm{\sigma}(t) = \bm{\gamma}(t) - \frac{\alpha(t)}{\overline{m}(t)} \bm{v}(t). 
$$
$NS_{\bm{\gamma}}[\bm{v}]$ is a cone (namely, the striction curve $\bm{\sigma}$ is constant) if and only if $\sigma(t) = 0$ for all $t \in I$, where 
$$
\sigma(t) = - \frac{d}{dt} \left( \frac{\alpha(t)}{\overline{m}(t)} \right). 
$$
By a straightforward calculation, we have
\begin{equation*}
NS_{\bm{\gamma}}[\bm{v}]_t (t, \lambda) \times NS_{\bm{\gamma}}[\bm{v}]_\lambda (t, \lambda) = \{ \alpha(t) + \lambda \overline{m}(t) \} \bm{w}(t), 
\end{equation*}
so that $\bm{\sigma}$ gives a parametrization of the singular value set of $NS_{\bm{\gamma}}[\bm{v}]$ under the assumption $\overline{m}(t) \neq 0$ for all $t \in I$. 
Furthermore, $(NS_{\bm{\gamma}}[\bm{v}], \bm{w}, \bm{v})$ is a framed surface, see Remark \ref{bishop_normal}. 
Since $\overline{\ell}(t) = 0$, we can take $\varphi(t, \lambda) = 0$ by Proposition \ref{proposition:normal_surfaces_condition}. 
Then we have the curvature $C^F = (J^F, K^F, H^F)$ of the framed surface $(NS_{\bm{\gamma}}[\bm{v}], \bm{w}, \bm{v})$, where 
$$
J^F(t, \lambda) = \alpha(t) + \lambda \overline{m}(t), \ K^F(t, \lambda) = 0, \ H^F(t, \lambda) = - \frac{1}{2} \overline{n}(t). 
$$
By Proposition \ref{pro:front_condition}, $NS_{\bm{\gamma}}[\bm{v}]$ is a front at a singular point $(t_0, \lambda_0) \in I \times \R$ if and only if $\alpha(t_0) + \lambda_0 \overline{m}(t_0) = 0$ and $\overline{n}(t_0) \neq 0$. 
\par 
We rewrite criteria that a singular point of $NS_{\bm{\gamma}}[\bm{v}]$ is a cuspidal edge, a swallowtail, and a cuspidal cross cap (cf. Theorem \ref{proposition:ce_sw_ccr}). 
\begin{corollary}\label{sing_dev}
Let $(\bm{\gamma}, \bm{\nu}_1, \bm{\nu}_2):I \to \R^3 \times \Delta$ be a framed curve with a Bishop frame $\{ \bm{v}, \bm{w}, \bm{\mu} \}$ and $\overline{m}(t) \neq 0$. 
Suppose that $(t_0, \lambda_0) \in S(NS_{\bm{\gamma}}[\bm{v}])$ is a non-degenerate singular point of $NS_{\bm{\gamma}}[\bm{v}]$. 
\begin{itemize}
\item[\rm (1)] Suppose that $\overline{n}(t_0) \neq 0$. 
Then we have the following. 
	\begin{itemize}
	\item[\rm (i)] The singular point $(t_0, \lambda_0)$ of $NS_{\bm{\gamma}}[\bm{v}]$ is a cuspidal edge if and only if $\sigma(t_0) \neq 0$. 
	\item[\rm (ii)] The singular point $(t_0, \lambda_0)$ of $NS_{\bm{\gamma}}[\bm{v}]$ is a swallowtail if and only if $\sigma(t_0) = 0$ and $\dot{\sigma}(t_0) \neq 0$. 
	\end{itemize}
	\item[\rm (2)] Suppose that $\overline{n}(t_0) = 0$. 
	Then the singular point $(t_0, \lambda_0)$ of $NS_{\bm{\gamma}}[\bm{v}]$ is a cuspidal cross cap if and only if $\dot{\overline{n}}(t_0) \neq 0$ and $\sigma(t_0) \neq 0$. 
\end{itemize}
\end{corollary}
\par 
We introduce families of functions on framed curves that are useful for the study of normal surfaces and their striction curves. 
For a framed curve $(\bm{\gamma}, \bm{\nu}_1, \bm{\nu}_2):I \to \R^3 \times \Delta$ with a Bishop frame $\{ \bm{v}, \bm{w}, \bm{\mu} \}$, we define a function $G[\bm{w}]:I \times \R^3 \to \R$ by
\[
G[\bm{w}](t, \bm{x}) = (\bm{x} - \bm{\gamma}(t)) \cdot \bm{w}(t). 
\]
We call $G[\bm{w}]$ a {\it support function} on $\bm{\gamma}$ with respect to $\bm{w}$. 
We denote that $g_{\bm{x}_0}[\bm{w}](t) = (\bm{x_0} - \bm{\gamma}(t)) \cdot \bm{w}(t)$ for any $\bm{x}_0 \in \R^3$. 
Then we have the following proposition. 
\begin{proposition}
Let $(\bm{\gamma}, \bm{\nu}_1, \bm{\nu}_2):I \to \R^3 \times \Delta$ be a framed curve with a Bishop frame $\{ \bm{v}, \bm{w}, \bm{\mu} \}$. 
Suppose that $\overline{m}(t_0) \neq 0$ and $\overline{n}(t_0) \neq 0$. 
Then we have the following. 
\begin{itemize}
\item[\rm (1)] $g_{\bm{x}_0}[\bm{w}](t_0) = 0$ if and only if there exist $\lambda, \mu \in \R$ such that $\bm{x}_0 - \bm{\gamma}(t_0) = \lambda \bm{v}(t_0) + \mu \bm{\mu}(t_0)$. 
\item[\rm (2)] $g_{\bm{x}_0}[\bm{w}](t_0) = \dot{g}_{\bm{x}_0}[\bm{w}](t_0) = 0$ if and only if there exists $\lambda \in \R$ such that $\bm{x}_0 - \bm{\gamma}(t_0) = \lambda \bm{v}(t_0)$. 
\item[\rm (3)] $g_{\bm{x}_0}[\bm{w}](t_0) = \dot{g}_{\bm{x}_0}[\bm{w}](t_0) = \ddot{g}_{\bm{x}_0}[\bm{w}](t_0) = 0$ if and only if 
$$
\bm{x}_0 - \bm{\gamma}(t_0) = - \frac{\alpha(t_0)}{\overline{m}(t_0)} \bm{v}(t_0). 
$$
\end{itemize}
\end{proposition}
\demo
Since $g_{\bm{x}_0}[\bm{w}](t) = (\bm{x_0} - \bm{\gamma}(t)) \cdot \bm{w}(t)$, we have the following calculations: 
\begin{align*}
{\rm (i)} \ g_{\bm{x}_0}[\bm{w}](t) &= (\bm{x_0} - \bm{\gamma}(t)) \cdot \bm{w}(t), \\
{\rm (ii)} \ \dot{g}_{\bm{x}_0}[\bm{w}](t) &= (\bm{x_0} - \bm{\gamma}(t)) \cdot ( \overline{n}(t) \bm{\mu}(t)), \\
{\rm (iii)} \ \ddot{g}_{\bm{x}_0}[\bm{w}](t) &= - \alpha(t) \overline{n}(t) + (\bm{x}_0 - \bm{\gamma}(t)) \cdot ( - \overline{m}(t) \overline{n}(t) \bm{v}(t) - \overline{n}^2(t) \bm{w}(t) + \dot{\overline{n}}(t) \bm{\mu}(t)  ). 
\end{align*}
By solving the above system of equations, we have assertions (1), (2), and (3). 
\enD
%
%
%
\section{Circular evolutes and involutes}\label{sec_evo_invo}
We introduce circular evolutes with respect to Bishop directions and involutes of framed curves. 
First, we introduce circular evolutes of framed curves with respect to Bishop directions. 
Let $(\bm{\gamma}, \bm{\nu}_1, \bm{\nu}_2):I \to \R^3 \times \Delta$ be a framed curve with a Bishop frame $\{ \bm{v}, \bm{w}, \bm{\mu} \}$, namely, $\overline{\ell}(t) = 0$ for all $t \in I$ in the Frenet-Serret type formula in \eqref{curvature-rotated}. 
We assume that $\overline{m}(t) \neq 0$ for all $t \in I$. 
Then we define a space curve $E_{\bm{\gamma}}[\bm{v}]:I \to \R^3$ by 
$$
E_{\bm{\gamma}}[\bm{v}](t) = \bm{\gamma}(t) - \frac{\alpha(t)}{\overline{m}(t)} \bm{v}(t). 
$$
This is nothing but the striction curve of $NS_{\bm{\gamma}}[\bm{v}]$, see section \ref{section:normal_developable_surfaces}. 
We call the space curve $E_{\bm{\gamma}}[\bm{v}]$ a ({\it circular}) {\it evolute} of $\bm{\gamma}$ with respect to $\bm{v}$ (or, {\it $\bm{v}$-evolute}). 
By a straightforward calculation, we have 
$$
\dot{E}_{\bm{\gamma}}[\bm{v}](t) = \dot{\bm{\gamma}}(t) - \frac{d}{dt} \left( \frac{\alpha(t)}{\overline{m}(t)} \right) \bm{v}(t) - \frac{\alpha(t)}{\overline{m}(t)} \dot{\bm{v}}(t) = - \frac{d}{dt} \left( \frac{\alpha(t)}{\overline{m}(t)} \right) \bm{v}(t),  
$$
so that $\dot{E}_{\bm{\gamma}}[\bm{v}](t) \cdot \bm{w}(t) = \dot{E}_{\bm{\gamma}}[\bm{v}](t) \cdot \bm{\mu}(t) = 0$ for all $t \in I$. 
Therefore, $(E_{\bm{\gamma}}[\bm{v}], \bm{w}, \bm{\mu}):I \to \R^3 \times \Delta$ is a framed curve with the curvature $(\overline{n},0,-\overline{m},-(d/dt)(\alpha/\overline{m}))$. 
Note that $\{ \bm{w}, \bm{\mu}, \bm{v} \}$ may not be a Bishop frame along $E_{\bm{\gamma}}[\bm{v}]$. 
We now take a Bishop frame along $ E_{\bm{\gamma}}[\bm{v}]$ by a rotated frame. 
We define vectors $\bm{v}_E$ and $\bm{w}_E$ by $\bm{v}_E(t) = \cos \theta_E(t) \bm{w}(t) - \sin \theta_E(t) \bm{\mu}(t)$ and $\bm{w}_E(t) = \sin \theta_E(t) \bm{w}(t) + \cos \theta_E(t) \bm{\mu}(t)$, where $\theta_E$ satisfies $\dot{\theta}_E(t) = \overline{n}(t)$. 
Then $(E_{\bm{\gamma}}[\bm{v}], \bm{v}_E, \bm{w}_E):I \to \R^3 \times \Delta$ is also a framed curve. 
Since $\bm{\mu}_E(t) = \bm{v}_E(t) \times \bm{w}_E(t) = \bm{v}(t)$, the curvature $(\ell_E, m_E, n_E, \alpha_E)$ of the framed curve $(E_{\bm{\gamma}}[\bm{v}], \bm{v}_E, \bm{w}_E)$ is given by 
\begin{eqnarray*}
\ell_E(t) &=& \dot{\bm{v}}_E(t) \cdot \bm{w}_E(t) = 0, \\
m_E(t) &=& \dot{\bm{v}}_E(t) \cdot \bm{\mu}_E(t) = \overline{m}(t) \sin \theta_E(t),\\ n_E(t) &=& \dot{\bm{w}}_E(t) \cdot 
\bm{\mu}_E(t) = - \overline{m}(t) \cos \theta_E(t), \\
\alpha_E(t) &=& \dot{E}_{\bm{\gamma}}[\bm{v}](t) \cdot \bm{\mu}_E(t) = - \frac{d}{dt} \left( \frac{\alpha(t)}{\overline{m}(t)} \right) = \sigma(t). 
\end{eqnarray*}
The above Bishop frame $\{ \bm{v}_E, \bm{w}_E, \bm{v} \}$ along $E_{\bm{\gamma}}[\bm{v}]$ is necessary to consider second circular evolutes, namely, circular evolutes of the circular evolute. 
\par 
We have the following relations between circular evolutes and parallel curves. 
\begin{proposition}
The $\bm{v}$-evolute of the $\bm{v}$-parallel curve of $\bm{\gamma}$ coincides with the $\bm{v}$-evolute of $\bm{\gamma}$, namely, $E_{P_{\bm{\gamma}}[\bm{v}]}[\bm{v}](t) = E_{\bm{\gamma}}[\bm{v}](t)$. 
\end{proposition}
\demo
By a straightforward calculation, we have 
\begin{align*}
E_{P_{\bm{\gamma}}[\bm{v}]}[\bm{v}](t) &= P_{\bm{\gamma}}[\bm{v}](t) - \frac{\alpha_P(t)}{m_P(t)} \bm{v}(t) \\
&= \bm{\gamma}(t) + \lambda \bm{v}(t) - \frac{\alpha(t) + \lambda \overline{m}(t)}{\overline{m}(t)} \bm{v}(t) \nonumber \\
&= \bm{\gamma}(t) - \frac{\alpha(t)}{\overline{m}(t)} \bm{v}(t) \nonumber \\
&= E_{\bm{\gamma}}[\bm{v}](t), \nonumber 
\end{align*}
so that the assertion holds. 
\enD
Second, we introduce involutes of framed curves. 
Let $(\bm{\gamma}, \bm{\nu}_1, \bm{\nu}_2):I \to \R^3 \times \Delta$ be a framed curve with $m^2(t) + n^2(t) \neq 0$ for all $t \in I$. 
Then we define a space curve $I_{\bm{\gamma}}[t_0]:I \to \R^3$ by 
$$
I_{\bm{\gamma}}[t_0](t) = \bm{\gamma}(t)  - \left( \int_{t_0}^t \alpha(t) \ dt  \right) \bm{\mu}(t) 
$$
for a fixed $t_0 \in I$. 
We call the space curve $I_{\bm{\gamma}}[t_0]$ an {\it involute} of $\bm{\gamma}$ with respect to $t_0$ (or, {\it $t_0$-involute}). 
We define smooth maps $\bm{\xi}$ and $\bm{\eta}:I \to S^2$ by 
$$
\bm{\xi}(t) = \frac{n(t) \bm{\nu}_1(t) - m(t) \bm{\nu}_2(t)}{\sqrt{m^2(t) + n^2(t)}}, \quad\bm{\eta}(t) = \bm{\xi}(t) \times \bm{\mu}(t) = \frac{-m(t) \bm{\nu}_1(t) - n(t) \bm{\nu}_2(t)}{\sqrt{m^2(t) + n^2(t)}}. 
$$
By a straightforward calculation, we have 
\begin{align*}
\dot{I}_{\bm{\gamma}}[t_0](t) &= \dot{\bm{\gamma}}(t) - \frac{d}{dt} \left( \int_{t_0}^t \alpha(t) \ dt \right) \bm{\mu}(t) - \left( \int_{t_0}^t \alpha(t) \ dt \right) \dot{\bm{\mu}}(t) \\
&= \left( \int_{t_0}^t \alpha(t) \ dt \right) \left( m(t) \bm{\nu}_1(t) + n(t) \bm{\nu}_2(t) \right), 
\end{align*}
so that $\dot{I}_{\bm{\gamma}}[t_0](t) \cdot \bm{\xi}(t) = \dot{I}_{\bm{\gamma}}[t_0](t) \cdot \bm{\mu}(t) = 0$ for all $t \in I$. 
Therefore,  $(I_{\bm{\gamma}}[t_0], \bm{\xi}, \bm{\mu}):I \to \R^3 \times \Delta$ is a framed curve with the curvature $(0,f,\sqrt{m^2+n^2},-(\int^t_{t_0} \alpha(t)dt))\sqrt{m^2+n^2})$, where 
$$
f(t)=\frac{\dot{m}(t)n(t)-m(t)\dot{n}(t)-(m^2(t)+n^2(t))\ell(t)}{m^2(t)+n^2(t)}.
$$ 
Note that $\{ \bm{\xi}, \bm{\mu}, \bm{\eta} \}$ is one of the Bishop frames along $I_{\bm{\gamma}}[t_0]$. 
We now take other Bishop frame along $I_{\bm{\gamma}}[t_0]$ by a rotated frame. 
We define vectors $\bm{v}_I$ and $\bm{w}_I$ by $\bm{v}_I(t) = \cos \theta_I \bm{\xi}(t) - \sin \theta_I \bm{\mu}(t)$ and $\bm{w}_I(t) = \sin \theta_I \bm{\xi}(t) + \cos \theta_I \bm{\mu}(t)$, where $\theta_I$ is a constant. 
Then $(I_{\bm{\gamma}}[t_0], \bm{v}_I, \bm{w}_I):I \to \R^3 \times \Delta$ is also a framed curve. 
Since $\bm{\mu}_I(t) = \bm{v}_I(t) \times \bm{w}_I(t) = \bm{\eta}(t)$, the curvature $(\ell_I, m_I, n_I, \alpha_I)$ of the framed curve $(I_{\bm{\gamma}}[t_0], \bm{v}_I, \bm{w}_I)$ is given by 
\begin{align*}
\ell_I(t) &= \dot{\bm{v}}_I(t) \cdot \bm{w}_I(t) = 0, \\
m_I(t) &= \dot{\bm{v}}_I(t) \cdot \bm{\mu}_I(t) = f(t) \cos \theta_I - \sqrt{m^2(t) + n^2(t)} \sin \theta_I, \\
n_I(t) &= \dot{\bm{w}}_I(t) \cdot \bm{\mu}_I(t)  = f(t) \sin \theta_I + \sqrt{m^2(t) + n^2(t)} \cos \theta_I, \\
\alpha_I(t) &= \dot{I}_{\bm{\gamma}}[t_0](t) \cdot \bm{\mu}_I(t) = - \left( \int_{t_0}^t \alpha(t) \ dt \right) \sqrt{m^2(t) + n^2(t)}. 
\end{align*}
The above Bishop frame $\{ \bm{v}_I, \bm{w}_I, \bm{\eta} \}$ along $I_{\bm{\gamma}}[t_0]$ is necessary to consider circular evolutes of the involute. 
\par
\begin{figure}[htbp]
\begin{minipage}{0.5\hsize}
\begin{center}
\includegraphics[width=50mm]{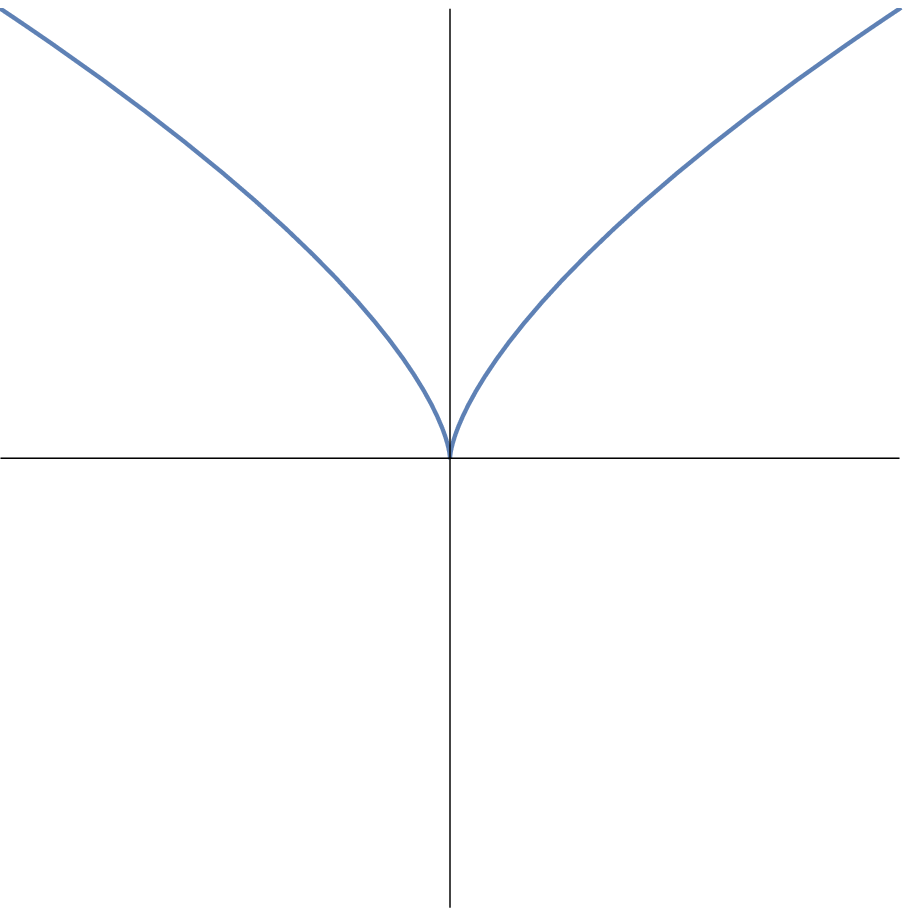}
\end{center}
\caption{A $3/2$-cusp $t \mapsto (t^2, t^3,0)$.}
\label{fig:2_3-cusp}
\end{minipage}
\begin{minipage}{0.5\hsize}
\begin{center}
\includegraphics[width=50mm]{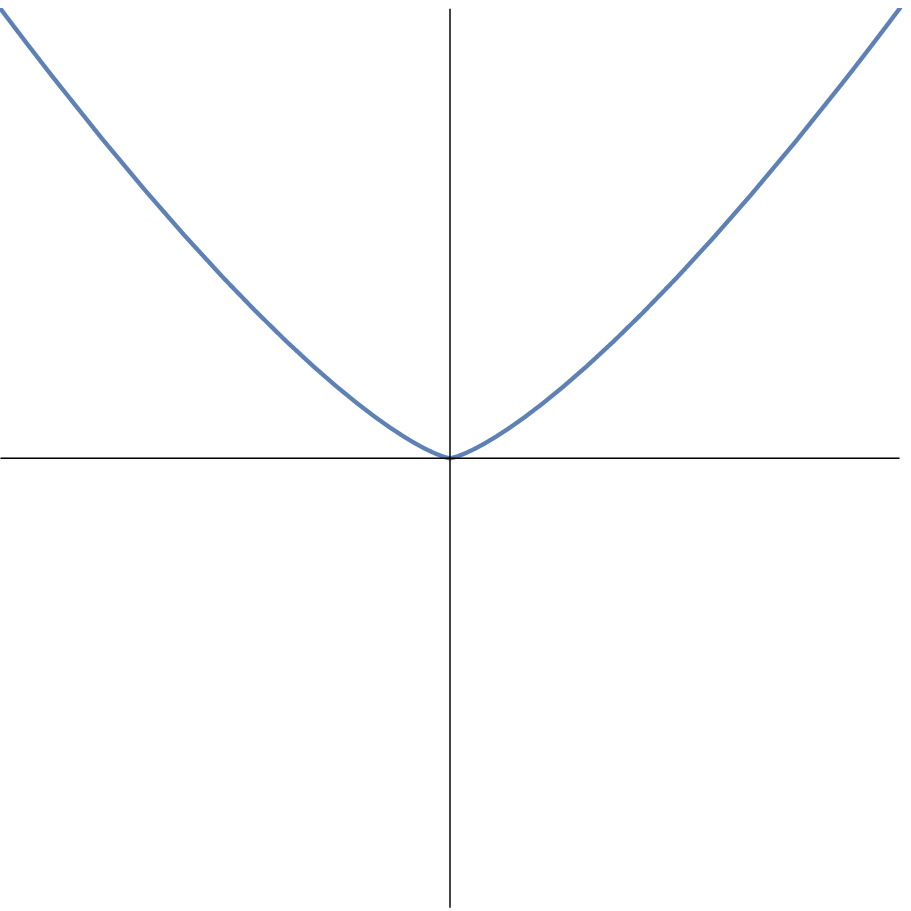}
\end{center}
\caption{A $4/3$-cusp $t \mapsto (t^3, t^4,0)$.}
\label{fig:3_4-cusp}
\end{minipage}
\end{figure}

 
Circular evolutes and involutes have some properties similar to plane curves, see \cite{FT3, FT4, FT5}. 
For instance, We have the following relations between circular evolutes and involutes. 
\begin{proposition}
Let $(\bm{\gamma}, \bm{\nu}_1, \bm{\nu}_2):I \to \R^3 \times \Delta$ be a framed curve with a Bishop frame $\{ \bm{v}, \bm{w}, \bm{\mu} \}$ and $\overline{m}(t) \neq 0$ for all $t \in I$. 
Then $I_{E_{\bm{\gamma}}[\bm{v}]}[t_0] = \bm{\gamma}(t) - (\alpha(t_0)/\overline{m}(t_0)) \bm{v}(t)=P_{\gamma}[\bm{v}](t)$, where $\lambda=- \alpha(t_0)/\overline{m}(t_0)$. 
\end{proposition}
\demo
By a straightforward calculation, we have 
\begin{align*}
I_{E_{\bm{\gamma}}[\bm{v}]}[t_0](t) &= E_{\bm{\gamma}}[\bm{v}](t) - \left( \int_{t_0}^t \alpha_E(t) \ dt  \right) \bm{v}(t) \\
&= \bm{\gamma}(t) - \frac{\alpha(t)}{\overline{m}(t)} \bm{v}(t) + \left\{ \int_{t_0}^{t}  \frac{d}{dt} \left( \frac{\alpha(t)}{\overline{m}(t)} \right) dt \right\} \bm{v}(t) \\
&= \bm{\gamma}(t) - \frac{\alpha(t)}{\overline{m}(t)} \bm{v}(t) + \left( \frac{\alpha(t)}{\overline{m}(t)} - \frac{\alpha(t_0)}{\overline{m}(t_0)} \right) \bm{v}(t) \\
&= \bm{\gamma}(t) -\frac{\alpha(t_0)}{\overline{m}(t_0)} \bm{v}(t), 
\end{align*}
so that the assertion holds. 
\enD
\begin{proposition}
Let $(\bm{\gamma}, \bm{\nu}_1, \bm{\nu}_2):I \to \R^3 \times \Delta$ be a framed curve with $m^2(t) + n^2(t) \neq 0$ for all $t \in I$. 
If we take $\theta_I = \pi/2$, then $E_{I_{\bm{\gamma}}[t_0]}[\bm{v}_I](t) = \bm{\gamma}(t)$ for any fixed $t_0 \in I$. 
\end{proposition}
\demo
If we take $\theta_I =\pi/2$, then $\bm{v}_I(t) = - \bm{\mu}(t) \quad \text{and} \quad m_I(t) = - \sqrt{m^2(t) + n^2(t)}$. 
By a straightforward calculation, we have
\begin{align*}
E_{I_{\bm{\gamma}}[t_0]}[\bm{v}_I](t) &= I_{\bm{\gamma}}[t_0](t) - \frac{\alpha_I(t)}{m_I(t)} \bm{v}_I(t) \\
&= \bm{\gamma}(t)  - \left( \int_{t_0}^t \alpha(t) \ dt  \right) \bm{\mu}(t) + \left( \int_{t_0}^t \alpha(t) \ dt  \right) \bm{\mu}(t) \\
&= \bm{\gamma}(t), 
\end{align*}
so that assertion holds. 
\enD
We now consider singularities of circular evolutes and involutes of framed curves. 
A singular point $t_0$ of a map $\bm{\gamma}:\R \to \R^3$ is called a {\it $3/2$-cusp}  (respectively, {\it $4/3$-cusp}) if the map-germ $\bm{\gamma}$ at $p$ is $\mathcal{A}$-equivalent to $t \mapsto (t^2, t^3,0)$ (respectively, $t \mapsto (t^3, t^4,0)$) at $0$, see Figures \ref{fig:2_3-cusp} and \ref{fig:3_4-cusp}. 
It is well-known that a singular point $t_0$ of $\bm{\gamma}$ is a $3/2$-cusp if and only if $\dot{\bm{\gamma}}(t_0) = 0$ and $\rank(\ddot{\bm{\gamma}}(t_0), \bm{\gamma}^{(3)}(t_0)) = 2$. 
Furthermore, a singular point $t_0$ of $\bm{\gamma}$ is a $4/3$-cusp if and only if $\dot{\bm{\gamma}}(t_0) = \ddot{\bm{\gamma}}(t_0) = 0$ and $\rank(\bm{\gamma}^{(3)}(t_0), \bm{\gamma}^{(4)}(t_0)) = 2$. 
\par 
By using the above criteria, we give the following relations among singularities of the original curve $\bm{\gamma}$,  the circular evolute $E_{\bm{\gamma}}[\bm{v}]$, and the involute $I_{\bm{\gamma}}[t_0]$. 
\begin{proposition}\label{prop:sing_evo}
Let $(\bm{\gamma}, \bm{\nu}_1, \bm{\nu}_2):I \to \R^3 \times \Delta$ be a framed curve with a Bishop frame $\{ \bm{v}, \bm{w}, \bm{\mu} \}$ and $\overline{m}(t) \neq 0$ for all $t \in I$. 
Suppose that $t_0$ is a singular point of $\bm{\gamma}$, namely, $\alpha(t_0) = 0$. 
Then we have the following. 
\begin{itemize}
\item[\rm (1)] The singular point $t_0$ of $\bm{\gamma}$ is a  $3/2$-cusp if and only if the point $t_0$ of $E_{\bm{\gamma}}[\bm{v}]$ is regular. 
\item[\rm (2)] The singular point $t_0$ of $\bm{\gamma}$ is a $4/3$-cusp if and only if the singular point $t_0$ of $E_{\bm{\gamma}}[\bm{v}]$ is a $3/2$-cusp. 
\end{itemize}
\end{proposition}
\demo
(1) By a straightforward calculation, we have $\dot{\bm{\gamma}}(t_0) = 0$ and $\rank (\ddot{\bm{\gamma}}(t_0), \bm{\gamma}^{(3)}(t_0)) = 2$ if and only if $\alpha(t_0) = 0$ and $\dot{\alpha}(t_0) \neq 0$. 
On the other hand, we have $\dot{E}_{\bm{\gamma}}[\bm{v}](t_0) \neq 0$ if and only if $\alpha_E(t_0) \neq 0$. 
Since $\alpha(t_0)=0$ and 
\begin{equation}\label{eq:1}
\alpha_E(t) = - \frac{1}{\overline{m}^2(t)} \left( \dot{\alpha}(t) \overline{m}(t) - \alpha(t) \dot{\overline{m}}(t) \right),
\end{equation}
assertion (1) holds. 
\par 
(2) By a straightforward calculation, we have $\dot{\bm{\gamma}}(t_0) = \ddot{\bm{\gamma}}(t_0) = 0$ and $\rank ( \bm{\gamma}^{(3)}(t_0), \bm{\gamma}^{(4)}(t_0) ) = 2$ if and only if $\alpha(t_0) = \dot{\alpha}(t_0) = 0$, and $\ddot{\alpha}(t_0) \neq 0$. 
On the other hand, we have $\dot{E}_{\bm{\gamma}}[\bm{v}](t_0) = 0$ and $\rank (\ddot{E}_{\bm{\gamma}}[\bm{v}](t_0), E_{\bm{\gamma}}[\bm{v}]^{(3)}(t_0)) = 2$ if and only if $\alpha_E(t_0) = 0$ and $\dot{\alpha}_E(t_0) \neq 0$. 
Since equations $\alpha(t_0)=0$, \eqref{eq:1} and 
$$
\dot{\alpha}_E(t) = \frac{2 \dot{\overline{m}}(t)}{\overline{m}^3(t)} \left( \dot{\alpha}(t) \overline{m}(t) - \alpha(t) \dot{\overline{m}}(t) \right) - \frac{1}{\overline{m}^2(t)} \left( \ddot{\alpha}(t) \overline{m}(t) - \alpha(t) \ddot{\overline{m}}(t) \right),
$$
assertion (2) holds. 
\enD
\begin{proposition}\label{prop:sing_invo}
Let $(\bm{\gamma}, \bm{\nu}_1, \bm{\nu}_2):I \to \R^3 \times \Delta$ be a framed curve with $m^2(t) + n^2(t) \neq 0$ for all $t \in I$. 
Suppose that $t_1$ is a singular point of $I_{\bm{\gamma}}[t_0]$. 
Then we have the following. 
\begin{itemize}
\item[\rm (1)] The point $t_1$ of $\bm{\gamma}$ is regular if and only if the singular point $t_1$ of $I_{\bm{\gamma}}[t_0]$ is a $3/2$-cusp. 
\item[\rm (2)] The singular point $t_1$ of $\bm{\gamma}$ is a $3/2$-cusp if and only if the singular point $t_1$ of $I_{\bm{\gamma}}[t_0]$ is a $4/3$-cusp. 
\end{itemize}
\end{proposition}
\demo
(1) By a straightforward calculation, we have $\dot{\bm{\gamma}}(t_1) \neq 0$ if and only if $\alpha(t_1) \neq 0$. 
On the other hand, we have $\dot{I}_{\bm{\gamma}}[t_0](t_1) = 0$ and $\rank (\ddot{I}_{\bm{\gamma}}[t_0](t_1),I_{\bm{\gamma}}[t_0]^{(3)}(t_1)  ) = 2$ if and only if $\alpha_I(t_1) = 0$ and $\dot{\alpha}_I(t_1) \neq 0$. 
Since
\begin{equation}\label{eq2}
\alpha_I(t) = - \left( \int_{t_0}^{t} \alpha(t) \ dt \right) \sqrt{m^2(t) + n^2(t)}
\end{equation}
and
\begin{equation}\label{eq3}
\dot{\alpha}_I(t) = - \alpha(t) \sqrt{m^2(t) + n^2(t)} - \left( \int_{t_0}^{t} \alpha(t) \ dt \right) \frac{d}{dt} \left( \sqrt{m^2(t) + n^2(t)} \right),
\end{equation}
assertion (1) holds. 
\par 
(2) By a straightforward calculation, we have $\dot{\bm{\gamma}}(t_1) = 0$ and $\rank (\ddot{\bm{\gamma}}(t_1), \bm{\gamma}^{(3)}(t_1)) = 2$ if and only if $\alpha(t_1) = 0$ and $\dot{\alpha}(t_1) \neq 0$. 
On the other hand, we have $\dot{I}_{\bm{\gamma}}[t_0](t_1) = \ddot{I}_{\bm{\gamma}}[t_0](t_1) = 0$ and $\rank ( I_{\bm{\gamma}}[t_0]^{(3)}(t_1), I_{\bm{\gamma}}[t_0]^{(4)}(t_1) ) = 2$ if and only if $\alpha_I(t_1) = \dot{\alpha}_I(t_1) = 0$ and $\ddot{\alpha}_I(t_1) \neq 0$. 
Since equations (\ref{eq2}), (\ref{eq3}) and 
\begin{align*}
\ddot{\alpha}_I(t) &= - \dot{\alpha}(t) \sqrt{m^2(t) + n^2(t)} - 2 \alpha(t) \frac{d}{dt} \left( \sqrt{m^2(t) + n^2(t)} \right) \\
& \qquad - \left( \int_{t_0}^{t} \alpha(t) \ dt \right) \frac{d^2}{dt^2} \left( \sqrt{m^2(t) + n^2(t)} \right), 
\end{align*}
assertion (2) holds. 
\enD
\par 
We now consider relations among normal surfaces, circular evolutes, and involutes of framed curves. 
Let $(\bm{\gamma}, \bm{\nu}_1, \bm{\nu}_2):I \to \R^3 \times \Delta$ be a framed curve with a Bishop frame $\{ \bm{v}, \bm{w}, \bm{\mu} \}$, namely, $\overline{\ell}(t) = 0$ for all $t \in I$ in the Frenet-Serret type formula in \eqref{curvature-rotated}. 
We assume that $\overline{m}(t) \neq 0$ for all $t \in I$. 
Then the $\bm{v}_E$-normal surface of $E_{\bm{\gamma}}[\bm{v}]$ is given by
$$
NS_{E_{\bm{\gamma}}[\bm{v}]}[\bm{v}_E](t, \lambda) = E_{\bm{\gamma}}[\bm{v}](t) + \lambda \bm{v}_E(t). 
$$
By results in section \ref{section:normal_developable_surfaces}, the singular value set of $NS_{E_{\bm{\gamma}}[\bm{v}]}[\bm{v}_E]$ is parametrized by the secondly circular evolute $E_{E_{\gamma}[\bm{v}]}[\bm{v}_E]$. 
\par 
On the other hand, let $(\bm{\gamma}, \bm{\nu}_1, \bm{\nu}_2):I \to \R^3 \times \Delta$ be a framed curve with $m^2(t) + n^2(t) \neq 0$ for all $t \in I$. 
If we take $\theta_I = \pi/2$ (namely, $\bm{v}_I(t) = - \bm{\mu}(t)$), then the $\bm{v}_I$-normal surface of $I_{\bm{\gamma}}[t_0]$ is given by 
$$
NS_{I_{\bm{\gamma}}[t_0]}[\bm{v}_I](t, \lambda) = I_{\bm{\gamma}}[t_0](t) + \lambda \bm{v}_I(t). 
$$
By results in section \ref{section:normal_developable_surfaces}, the singular value set of $NS_{I_{\bm{\gamma}}[t_0]}[\bm{v}_I]$ is parametrized by the original curve $\bm{\gamma}$. 
Since $NS_{I_{\bm{\gamma}}[t_0]}[\bm{v}_I]$ is a non-cylindrical developable surface, $NS_{I_{\bm{\gamma}}[t_0]}[\bm{v}_I](I \times \R)$ is the tangent developable surface of $\bm{\gamma}$ (cf. \cite{IT2}). 
Singularities of tangent developable surfaces are investigated in \cite{I1, I2, I3}. 
By the above properties, we have the diagram of Figure \ref{fig:diagram} (dotted arrows establish when we consider involutes with respect to singular points): 
\begin{figure}[htbp]
\begin{center}
\includegraphics[width=0.8\linewidth]{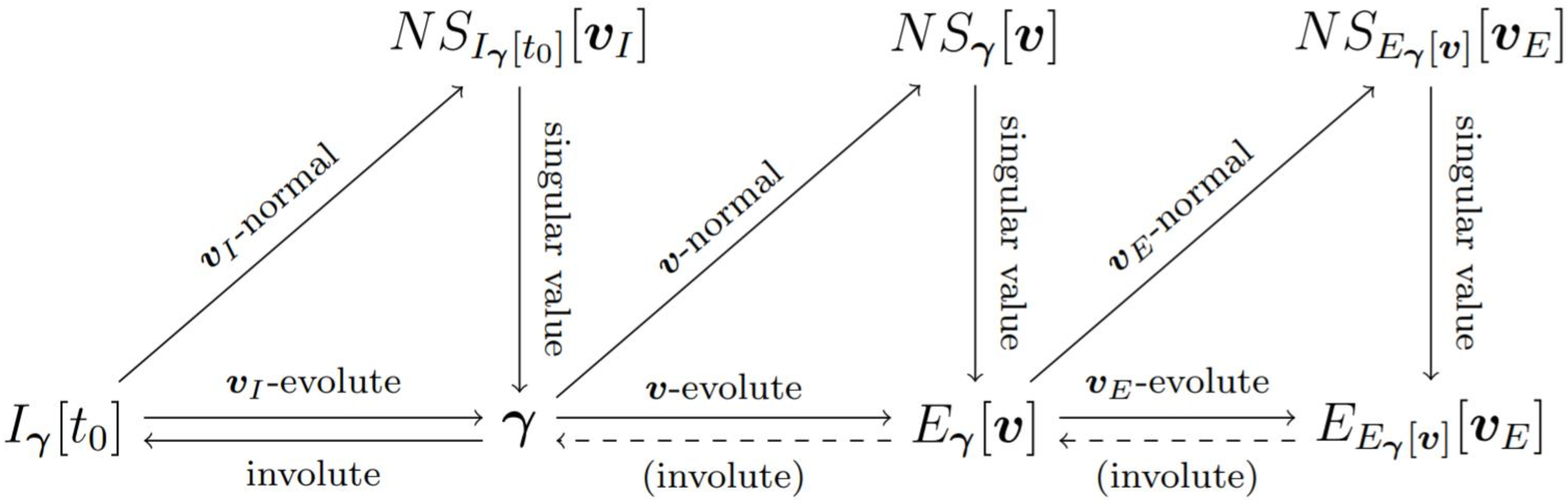}
\end{center}
\caption{Diagram of relations among normal surfaces, circular evolutes, and involutes.}
\label{fig:diagram}
\end{figure}
\begin{example}\label{example:last}
\rm We consider the spherical nephroid $\bm{\gamma}$ in Example \ref{example:spherical_nephroid}. 
Let $(\bm{\gamma}, \bm{\nu}_1, \bm{\nu}_2):[0, 2\pi) \to S^2 \times \Delta$ be 
\begin{align*}
\bm{\gamma}(t) &= \left( \frac{3}{4} \cos t - \frac{1}{4} \cos 3t, \frac{3}{4} \sin t - \frac{1}{4} \sin 3t, \frac{\sqrt{3}}{2} \cos t \right), \\
\bm{\nu}_1(t) &= \left( - \frac{3}{4} \sin t - \frac{1}{4} \sin 3t, \cos^3 t, \frac{\sqrt{3}}{2} \sin t \right), \\
\bm{\nu}_2(t) &= \left( \frac{3}{4} \cos t - \frac{1}{4} \cos 3t, \sin^3 t, \frac{\sqrt{3}}{2} \cos t \right). 
\end{align*}
Then singular points $t = 0$ and $\pi$ of $\bm{\gamma}$ are 3/2-cusps. 
If we take a Bishop frame $\{ \bm{v}, \bm{w}, \bm{\mu} \}$ by $\theta(t) = 0$ (namely, $\bm{v}(t) = \bm{\nu}_1(t)$ and $\bm{w}(t) = \bm{\nu}_2(t)$) and $t_0 = 0$, then we have  
\begin{align*}
NS_{\bm{\gamma}}[\bm{v}](t, \lambda) = \bm{\gamma}(t) + \lambda \bm{v}(t), \quad E_{\bm{\gamma}}[\bm{v}](t) = \bm{\gamma}(t) + \tan t \bm{v}(t), 
\end{align*}
and
\[
I_{\bm{\gamma}}[0](t) = \bm{\gamma}(t) - \sqrt{3} \left( 1 - \cos t \right) \bm{\mu}(t).
\]
By corollary \ref{sing_dev}, propositions \ref{prop:sing_evo} and \ref{prop:sing_invo}, we have the following: 
\begin{enumerate}
\item The singular points $(t, \lambda) = (0, 0)$ and $(\pi, 0)$ of $NS_{\bm{\gamma}}[\bm{v}]$ are cuspidal cross caps. 
\item The singular points $(t, \lambda) = (0, 0)$ and $(\pi, 0)$ of $NS_{I_{\bm{\gamma}}[0]}[\bm{\mu}]$ are swallowtails. 
\item The singular point $t = 0$ of $I_{\bm{\gamma}}[0]$ is a 4/3-cusp. 
\end{enumerate}
%
Moreover, by the above diagram, we can see relations among normal surfaces, circular evolutes, and involutes in Figures \ref{fig:np_evo} and \ref{fig:np_invo}. 
\begin{figure}[htbp]
\begin{minipage}{0.325\hsize}
\begin{center}
\includegraphics[width=40mm]{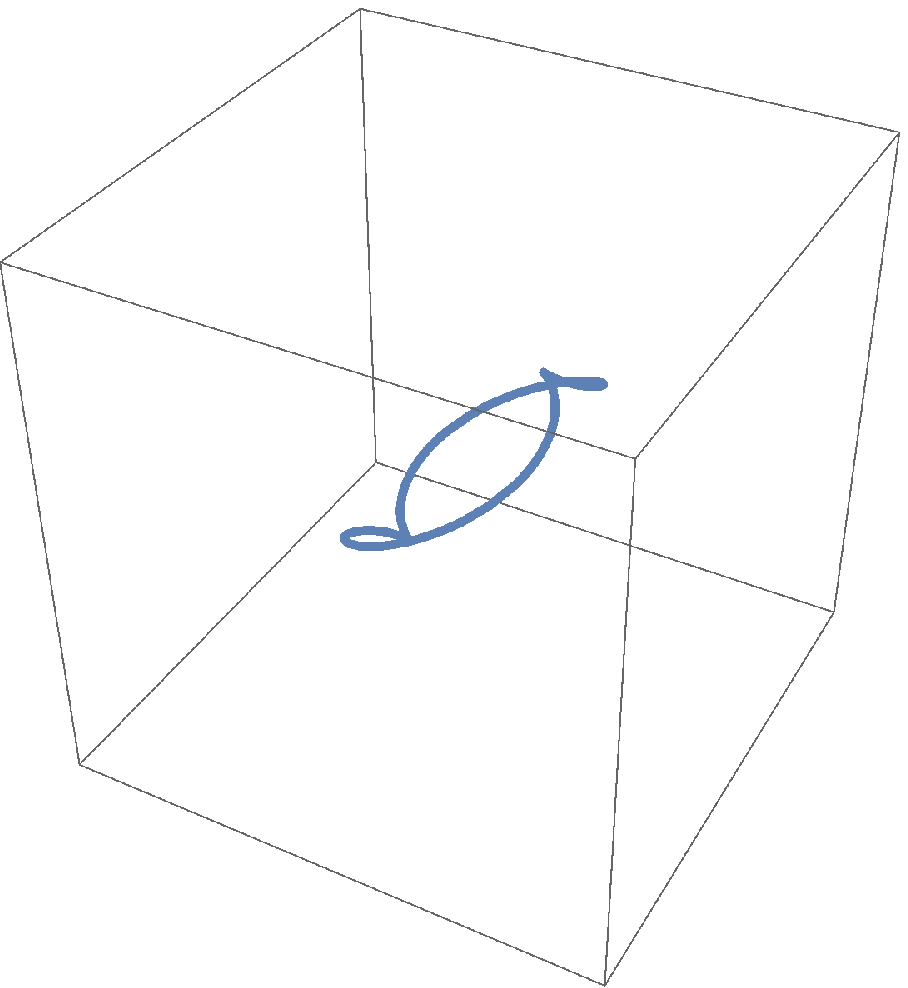}
\end{center}
\end{minipage}
\begin{minipage}{0.325\hsize}
\begin{center}
\includegraphics[width=40mm]{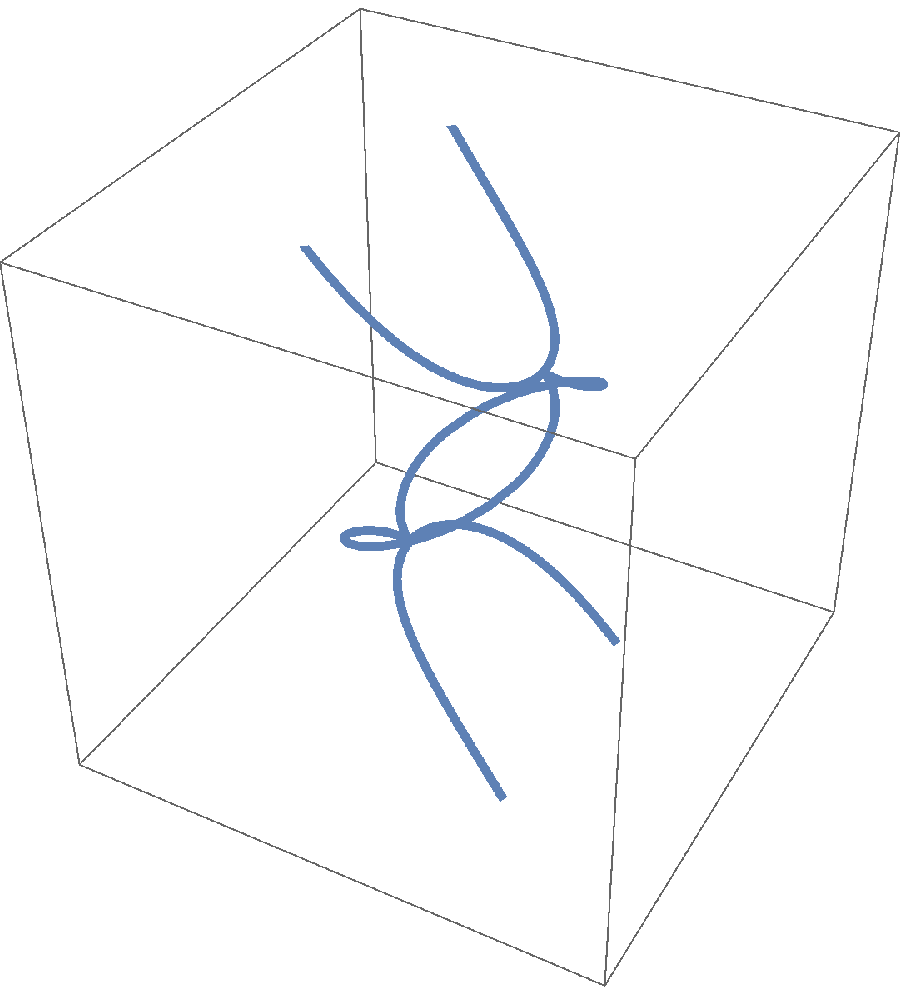}
\end{center}
\end{minipage}
\begin{minipage}{0.325\hsize}
\begin{center}
\includegraphics[width=40mm]{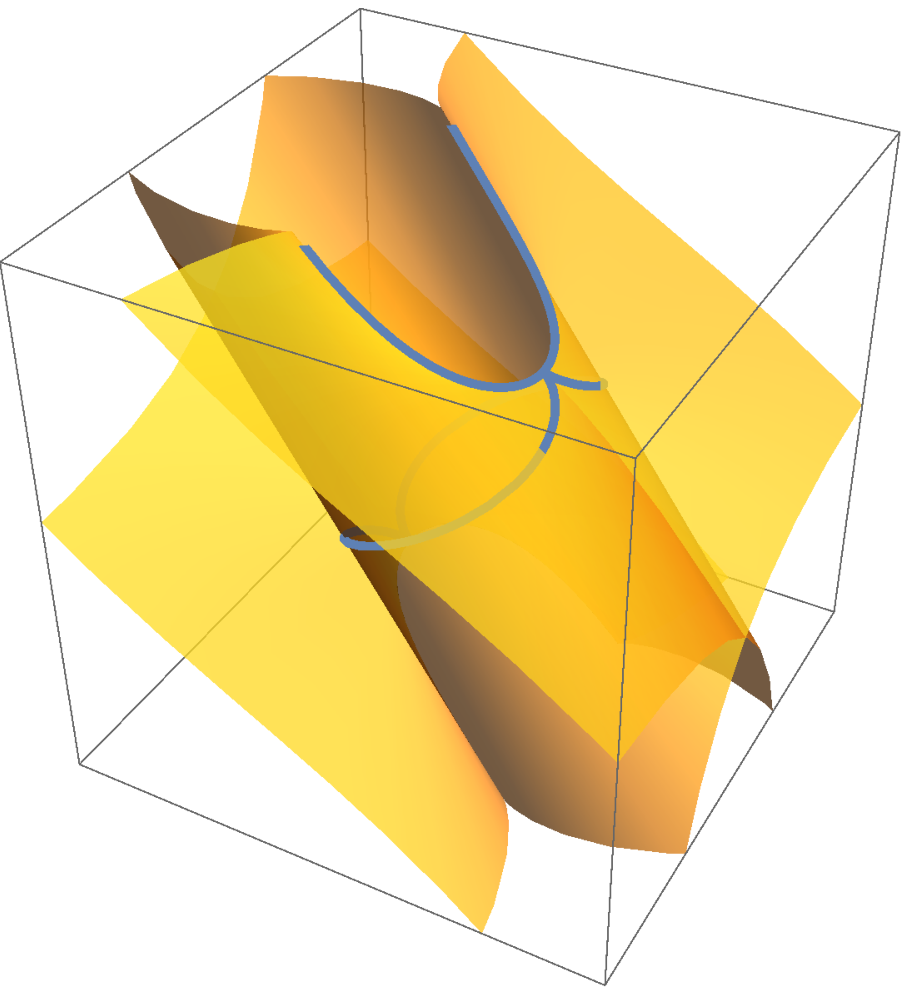}
\end{center}
\end{minipage}
\caption{Left to right: $\bm{\gamma}$, $(\bm{\gamma}, E_{\bm{\gamma}}[\bm{v}])$, and $(\bm{\gamma}, E_{\bm{\gamma}}[\bm{v}], NS_{\bm{\gamma}}[\bm{v}])$ in Example \ref{example:last}.}
\label{fig:np_evo}
\end{figure}
\begin{figure}[htbp]
\begin{minipage}{0.325\hsize}
\begin{center}
\includegraphics[width=40mm]{eps/np_gamma.eps}
\end{center}
\end{minipage}
\begin{minipage}{0.325\hsize}
\begin{center}
\includegraphics[width=40mm]{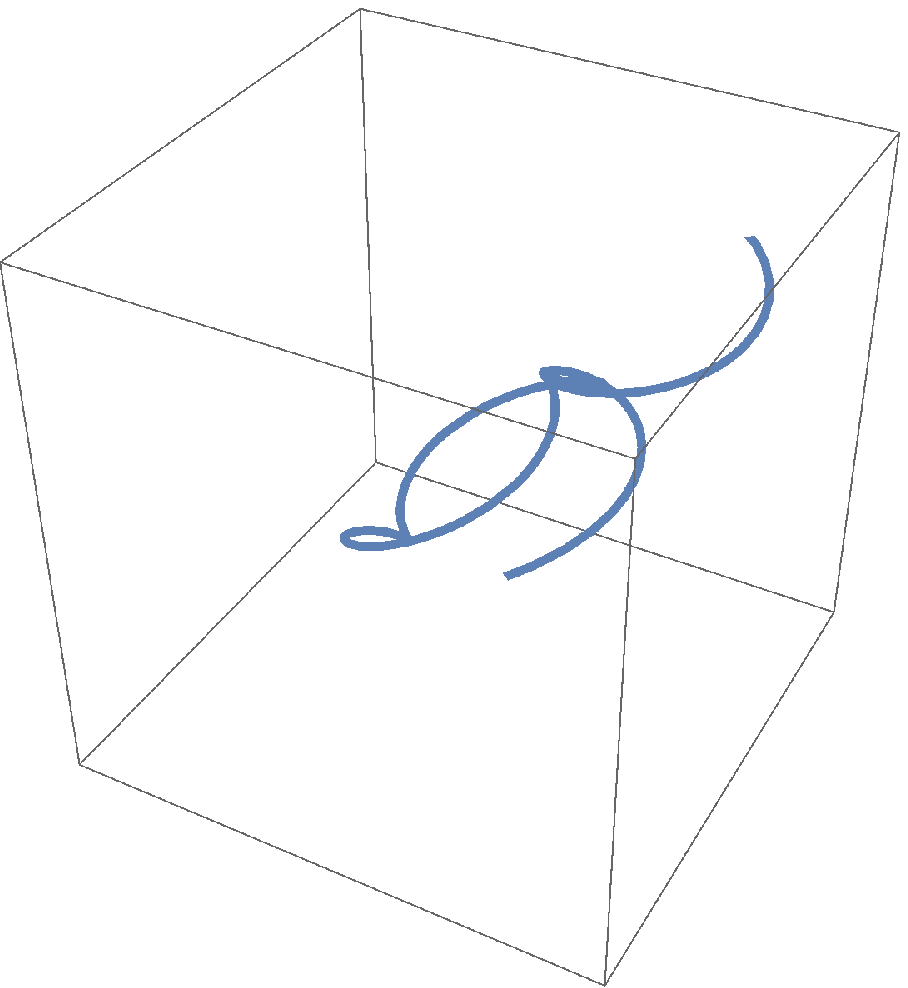}
\end{center}
\end{minipage}
\begin{minipage}{0.325\hsize}
\begin{center}
\includegraphics[width=40mm]{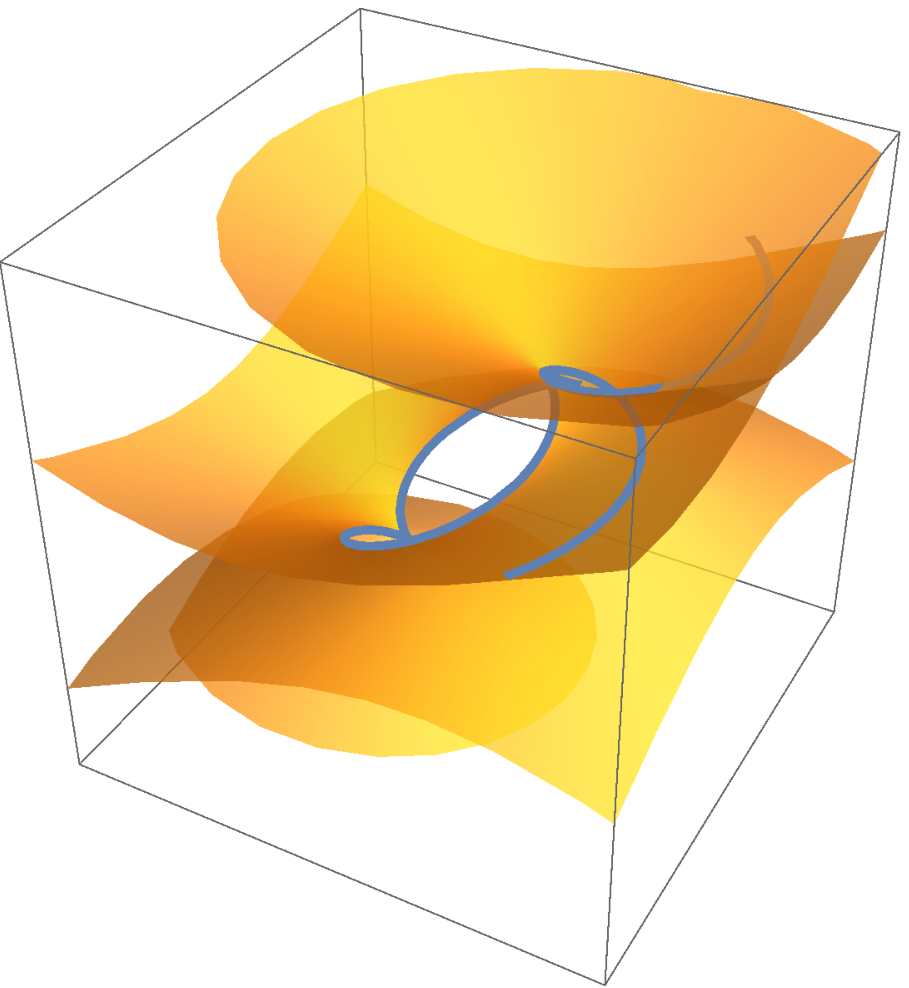}
\end{center}
\end{minipage}
\caption{Left to right: $\bm{\gamma}$, $(\bm{\gamma}, I_{\bm{\gamma}}[0])$, and $(\bm{\gamma}, I_{\bm{\gamma}}[0], NS_{I_{\bm{\gamma}}[0]}[\bm{v}_I])$ in Example \ref{example:last}.} 
\label{fig:np_invo}
\end{figure}
\end{example}

{\small 

Shun'ichi Honda,
\par
Chitose Institute of Science and Technology, Chitose 066-8655, Japan,
\par
E-mail address: s-honda@photon.chitose.ac.jp
\\

Masatomo Takahashi, 
\par
Muroran Institute of Technology, Muroran 050-8585, Japan,
\par
E-mail address: masatomo@mmm.muroran-it.ac.jp
}
\end{document}